\documentclass{elsarticle-modified}

\usepackage{lineno,hyperref}
\usepackage{amssymb,amsmath,amsfonts}
\usepackage{bm}
\usepackage[utf8]{inputenc}
\usepackage[T1]{fontenc}
\usepackage{verbatim}
\usepackage{xcolor}
\usepackage{mathrsfs}
\usepackage{cancel}
\usepackage{boxedminipage2e}
\modulolinenumbers[5]

\numberwithin{equation}{section}

\journal{J. Comput. Appl. Math.}

\newcommand{\x}[1]{{\color{red}{#1}}}
\renewcommand\t{\tau}
\newcommand\nl[1]{\\[#1\jot]}

\begin{document}

\newcommand\clas{{\bm c}}
\newcommand\symm{{\bm s}}

\newtheorem{remark}{Remark}
\newtheorem{theorem}{Theorem}
\newcommand\ol[1]{\overline{#1}}
\newcommand\nool[1]{#1}
\newcommand\ul[1]{\underline{#1}}
\newcommand\noul[1]{#1}
\newcommand\ub[2]{\underbrace{#1}_{#2}}
\newcommand\noub[2]{#1}
\newcommand\ca[1]{\cancel{#1}}
\newcommand\noca[2]{#1}
\newcommand{\CCC}{{{\mathbb{C}}\vphantom{|}}}
\newcommand{\NN}{{{\mathbb{N}}\vphantom{|}}}
\newcommand{\RR}{{{\mathbb{R}}\vphantom{|}}}
%
%
\newcommand{\dd}{{\mathrm d}}
\newcommand{\ee}{{\mathrm e}}
\newcommand{\ii}{{\mathrm i}}
\newcommand\Order{{\mathscr{O}}}
\newcommand\h{\frac12}
\renewcommand\th{\tfrac{1}{2}}
\newcommand\nEE[1]{\ee^{\,{#1}}}
\newcommand\nA{\mathcal{A}}
\newcommand\nB{\mathcal{B}}
\newcommand\nC{\mathcal{C}}
\newcommand\nE{\mathcal{E}}
\newcommand\nD{{\mathcal D}}
\newcommand\nF{{\mathcal F}}
\newcommand\nG{{\mathcal G}}
\newcommand\nL{{\mathcal L}}
\newcommand\nLL{{\tilde\nL}}
\newcommand\nS{{\mathcal S}}
\newcommand\nSD{{\widehat\nS}}
\newcommand\nLD{{\widehat\nL}}
\newcommand\T{{\mathcal T}}
\newcommand\TT{\widehat{\mathcal T}}
\newcommand\U{{\mathcal U}}
\newcommand\UU{\widehat{\mathcal U}}
\newcommand\V{{\mathcal V}}
\newcommand\nX{{\mathcal X}}
\newcommand\Y{{\mathcal Y}}
\newcommand\Z{{\mathcal Z}}

\newcommand\e[2]{\tt{#1}\,e{#2}}
\newcommand\f[1]{\tt{#1}}
\newcommand\void{~}
\newcommand\vv{$\vphantom{\int_x^x}$}
\newcommand\vvv{$\vphantom{\int_X^X}$}

\newcommand\Hfull{H_{\rm full}}
\newcommand\Hdiag{H_{\rm diag}}
\newcommand\Hsymm{H_{\rm symm}}
\newcommand\Hanti{H_{\rm anti}}
\def\C{\mathbb{C}}
\def\nsit{N}
\newcommand\dt{\tfrac{\dd}{\dd\t}}
\renewcommand{\t}{\tau}
\newcommand{\nP}{\mathcal{P}}

\newcommand\canc[1]{\text{$#1$}}
\newcommand\dtau{\tfrac{\dd}{\dd \tau}}
\newcommand\pdone{\partial_1}
\newcommand\Id{\text{Id}}
\newcommand\nDc{{\mathcal D}_{\clas}}
\newcommand\nDs{{\mathcal D}_{\symm}}
\newcommand\pdtwo{\partial_2}
\renewcommand\th{\tfrac{1}{2}}
\newcommand\nR{{\mathcal R}}
\newcommand\pdtnot{\tfrac{\partial}{\partial t_0}}
\newcommand\pdtau{\tfrac{\partial}{\partial \tau}}

\newcommand\rev[1]{{\color{red}{#1}}}
\newcommand\frage[1]{{\color{red}{#1}}}
\newcommand\checken[1]{{\color{green}{#1}}}
\newcommand\ak[1]{{\color{blue}{AK: #1}}}
\newcommand\kh[1]{{\color{blue}{KH: #1}}}

\begin{frontmatter}
\title{Efficient Magnus-type integrators for solar energy conversion in Hubbard models}
\author[tuwien1]{Winfried Auzinger}
\ead{winfried.auzinger@tuwien.ac.at}
\ead[url]{http://www.asc.tuwien.ac.at/~winfried/}
\author[tuwien1]{Juliette Dubois}
\ead{e11831465@student.tuwien.ac.at}
\author[tuwien2]{Karsten Held}
\ead{held@ifp.tuwien.ac.at}
\ead[url]{https://www.ifp.tuwien.ac.at/cms/}
\author[univie]{Harald Hofst{\"a}tter}
\ead{hofi@harald-hofstaetter.at}
\ead[url]{http://harald-hofstaetter.at}
\author[tuwien1]{Tobias Jawecki}
\ead{tobias.jawecki@tuwien.ac.at}
\author[tuwien2]{Anna Kauch}
\ead{kauch@ifp.tuwien.ac.at}
\ead[url]{https://www.ifp.tuwien.ac.at/cms/}
\author[univie]{Othmar Koch\corref{ourcorrespondingauthor}}
\ead{othmar@othmar-koch.org}
\ead[url]{http://www.othmar-koch.org}
\cortext[ourcorrespondingauthor]{Corresponding author}
\author[gdansk]{Karolina Kropielnicka}
\ead[url]{https://mat.ug.edu.pl/~kmalina/}
\ead{kmalina@mat.ug.edu.pl}
\author[bath]{Pranav Singh}
\ead{ps2106@bath.ac.uk}
\ead[url]{https://www.pranavsingh.co.uk/}
\author[tuwien2]{Clemens Watzenb\"ock}
\ead{clemens.watzenboeck@tuwien.ac.at}
\ead[url]{https://www.ifp.tuwien.ac.at/cms/}
\address[tuwien1]{Institute of Analysis and Scientific Computing, TU Wien, \\
               Wiedner Hauptstra{\ss}e~8-10, A--1040 Wien, Austria.}
\address[tuwien2]{Institute of Solid State Physics, TU Wien, \\
               Wiedner Hauptstra{\ss}e~8-10, A--1040 Wien, Austria.}
\address[univie]{Institut f{\"u}r Mathematik, Universit{\"a}t Wien, \\
                    Oskar-Morgensternplatz 1, A-1090 Wien, Austria.}
\address[gdansk]{Institute of Mathematics, University of Gdansk, ul. Wit Stwosz 57, 80-308 Gdansk}
\address[bath]{Department of Mathematical Sciences, University of Bath, Bath, BA2 7AY, U.K.}

\begin{abstract}
 Strongly interacting electrons in solids are generically described by Hubbard-type models, and the impact of solar light  can be modeled by an additional time-dependence. This
 yields a finite dimensional system of ordinary differential equations (ODE)s of Schr{\"o}dinger type, which can be solved numerically  by exponential
time integrators of Magnus type. The efficiency may be enhanced by combining
these with operator splittings. We will discuss several different
approaches of employing exponential-based methods in conjunction
with an adaptive Lanczos method for the evaluation of matrix exponentials
and compare their accuracy and efficiency.
For each integrator, we use defect-based local error estimators to
enable adaptive time-stepping.
This serves to reliably control the approximation error and reduce the computational effort.
\end{abstract}
\end{frontmatter}


\sloppy

\section{Introduction} \label{sec:intro}

The time evolution of a quantum mechanical system is generally described by a system of linear ordinary differential equations (ODE)'s
of Schr{\"o}dinger type
\begin{equation}\label{eq0}
\psi'(t) = - \mathrm{i} H(t) \psi(t) =: A(t) \psi(t), \qquad \psi(t_0)=\psi_0,
\end{equation}
with a large time-dependent Hermitian system matrix $H(t)$, state vector $\psi(t)$ at time $t$, and its derivative $\psi'(t)$.
The exact flow of \eqref{eq0} is denoted by $\psi(t)=\nE(t;t_0)\psi_0,$
and depends on the initial state $\psi_0$ and time $t_0$.
We will focus on the  movement and interaction of
electrons within Hubbard-type models, with the time dependence
originating from electric fields associated with a photon in the process of solar energy conversion~\cite{Werner2014,Innerberger2020}.
These models naturally have a discrete, and for a finite number of lattice sites finite, basis set in which the matrix $H(t)$ can be expressed.

The present study is motivated by a recent application where the efficient time propagation of models of
the type (\ref{eq0}) is of paramount importance: The simulation
of oxide solar cells with the goal of finding candidates for
new materials promising a gain in the solar cells' efficiency~\cite{kh:assmann13,Werner2014,Sorantin2018}.
For traditional materials such as silicon the efficiency of solar cells is fundamentally limited to 34\% due to the
\emph{Schockley--Queisser limit} \cite{SchockleyQueisser61}.
For overcoming this limit and building highly efficient solar cells,   \cite{kh:assmann13}
recently proposed a (for this purpose new) class of materials:
oxide heterostructures. These consist of at least two different
transition metal oxides. One of these can be the cheap and commonly used substrate material,
the perovskite SrTiO$_3$.
On top of this, one can stack layers of, e.g., LaVO$_3$ which has a preferable bandgap of $\sim 1.1\,$eV
for photovoltaic applications.
The equilibrium calculations of  \cite{kh:assmann13} have indicated the great potential of these oxide
heterostructures for solar cells;
and a LaVO$_3$/SrTiO$_3$ solar cell  was experimentally realized in \cite{Wang}, and likewise a LaFeO$_3$/SrTiO$_3$ one \cite{Nakamura2016}. A particular advantage of these oxide heterostructure solar cells is that one photon may excite two electron-hole pairs through a secondary process called impact ionization \cite{Manousakis2010,Werner2014,Sorantin2018}. This might serve to overcome the Schockley--Queisser limit. However,
for actually improving the efficiency of these solar cells, a better understanding
and calculation of the nonequilibrium processes, which are at the heart of the solar energy conversion,
are direly needed. Indeed, so far the efficiency of the solar cell has not been calculated as this results from an inherently nonequilibrium process.

In this paper, we compare numerical time integrators for the efficient approximation of the full dynamics of solar energy conversion in
 Hubbard models on a finite number of lattice sites.
The methods which are commonly considered as most appropriate for this task are based on the
matrix exponential function. When the latter is suitably approximated, the time propagation
conserves the norm of the wave function, which is the case for the Lanczos method we employ, see Section~\ref{sec:lanczos}.
This is violated, however, for popular one-step methods such as Runge--Kutta, which is
used for the purpose of comparisons in this study only. Our emphasis is on
adaptive time-stepping based on asymptotically correct estimates of the local error.
The advantage of adaptive step-size choice lies
not only in the potential for increased efficiency when the smoothness of
the solution varies over time, but more importantly in the reliable control
of the accuracy. The computational effort for the
evaluation of the error estimate may not always be compensated for by the optimal
choice of the local step-size. However, an optimal equidistant step-size
cannot usually be guessed a priori, while adaptive step-size selection
automatically adapts the time-steps such that the prescribed error tolerance
is satisfied. Hence, the accuracy of the numerical approximation
is reliably controlled.

The methods which we focus on are commutator-free Magnus-type methods \cite{alvfeh11,alvermanetal12},
but also classical Magnus integrators are considered \cite{haireretal02b,magnus54}.
Furthermore, a Magnus-Lanczos integrator employing operator splitting, inspired
by \cite{karo18a}, is tested. Explicit
Runge--Kutta methods are used for the purpose of comparisons.

In Section~\ref{sec:model} we describe the simple model of solar energy conversion in oxide solar cells
which we investigate. Section~\ref{sec:num-approach} describes the
numerical methods which we consider, and Section~\ref{sec:numexp} gives
the results of our comparative tests. \ref{sec:adapt} provides a
detailed description of the construction and implementation of
a new splitting-based Magnus--Strang integrator designed especially for
problems of the structure we are confronted with. Finally,
\ref{further} presents some supplementary numerical tests.

\section{The model} \label{sec:model}
For the description of the Hubbard model we resort to the \emph{second-quantization} formalism.
The Hubbard model first appears in~\cite{hubbard63,Gutzwiller63,Kanamori63} and since then became the basic model for describing strongly interacting (strongly correlated) electrons \cite{Ma93,PKBS16}.
It describes the electron occupation on a given number of sites, corresponding to Wannier discretization. Only a  single orbital per site is considered which allows for four states per site (no electron, one electron with spin-up or -down, two electrons). If there are two electrons on the same site this costs a Coulomb interaction $U$.  The Hubbard Hamiltonian for arbitrary hopping $v_{ij}$ reads
\begin{equation}\label{eq.HubHam}
H=\frac{1}{2} \sum_{i,j,\sigma} v_{ij} \hat{c}_{j\sigma}^{\dagger} \hat{c}_{i\sigma}^{\phantom\dagger}
+ \frac{1}{2}\sum_{i,\sigma} U \hat{n}_{i\sigma}\hat{n}_{i\bar{\sigma}},
\end{equation}
where $i,j$ sum over all $N$ sites  and the spins $\sigma,\sigma' \in\{\uparrow,\downarrow \}$ are either
\emph{up} or \emph{down}, and $\bar{\sigma}$ is the spin opposite to $\sigma$.
The notation $\hat{c}_{j\sigma}^\dagger c_{i\sigma}^{\phantom\dagger}$ describes a ``hopping'' from site $i$ to $j$ with  \emph{creation} and \emph{annihilation operators} $\hat{c}_{j \sigma}^\dagger$ and $\hat{c}_{i \sigma}$, respectively.
The \emph{hopping amplitudes} $v_{ij}$ with $i,j=1,\ldots,\nsit$ give the probability (rate) of such an electron hopping;
$\hat{n}_{j\sigma}=\hat{c}_{j\sigma}^\dagger \hat{c}_{j\sigma}^{\phantom\dagger}$ is the \emph{occupation number operator}, which counts the number of electrons with spin $\sigma$ at site $j$.
A derived observable which we will use later is the \emph{mean double occupation}
$\langle \hat d(t) \rangle = \frac1{N} \sum_{i=1}^{N} \langle \psi(t)|\hat d_i|\psi(t)\rangle,$
where the expectation value of $\hat d_i = \hat n_{i\uparrow}\hat n_{i \downarrow}$
is  1 (0) if there are  two (zero or one) electrons on site $i$.
For details on the notation in~\eqref{eq.HubHam} we recommend several references, e.g.~\cite{hubbard63,Ma93,PKBS16,Ja08}.

The time-dependence in the Hamiltonian (\ref{eq.HubHam}) is introduced through the photon
which excites the system out of equilibrium.
We approximate the photon by a classical electric field pulse of the form
$\vec{E}(t)=\vec{E}_0\ee^{-(t-t_p)^2/\sigma_p^2}\sin(\omega(t-t_p))$,
where $\omega$ denotes the frequency (energy) of the photon, $\sigma_p$ the width in time and $t_p$
the point in time of its impact.
We can relate  $\vec{E}(t)=-\partial_t \vec{A}(t)$ to a vector  potential $\vec{A}(t)$
using a  gauge without scalar potential. The vector potential in turn can be related, within the
Peierls' approximation~\cite{Peierls1933}, to a modified hopping amplitude \cite{Freericks06,Aoki2013,Werner2014}
\begin{equation}
v_{ij}\rightarrow v_{ij}(t) = v_{ij}\ee^{\mathrm{i} \int_{{\mathbf R}_i}^{{\mathbf R}_j}\vec{A}
({\mathbf r'},t)\,\dd{\mathbf r'}},
\label{eq:peierls}
\end{equation}
where ${\mathbf R}_i$ is the position of lattice site $i$.

\begin{figure}
	\centering
	\includegraphics[width=0.7\linewidth]{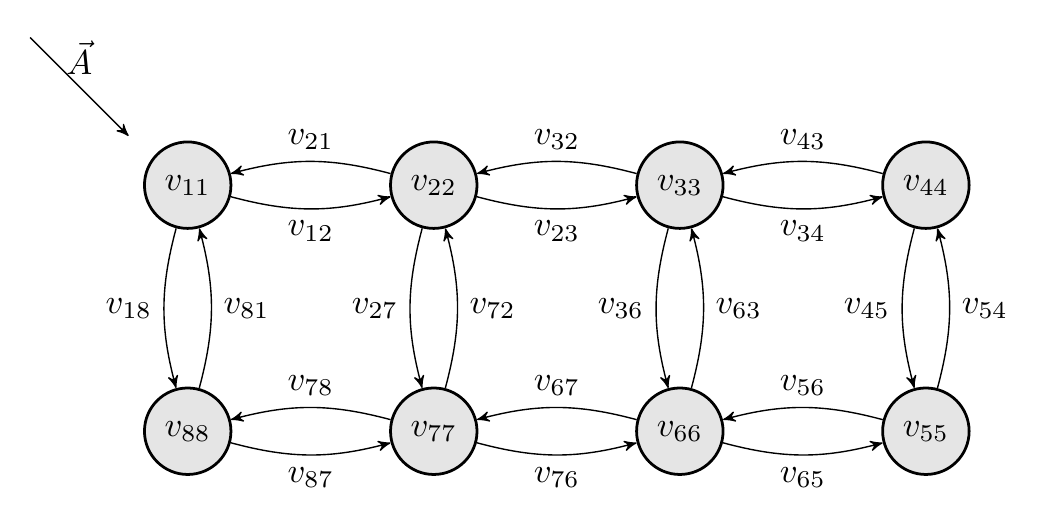}
	\caption{Geometry of a two-dimensional $8$-site lattice ($2\times 4$ ladder) with lexicographical ordering of the sites. The energies $v_{ii}$ describe the on-site potential, $v_{ij}$ describe the hoppings between sites $i$ and $j$. The vector potential $\vec{A}$ of the electromagnetic field is in the lattice plane and in the diagonal direction.}
	\label{fig:geometry}
\end{figure}

For our numerical tests we use Hubbard models with different geometric settings.
\begin{enumerate}
\item First, we model $8$ electrons on $8$ sites ($\nsit=8$) arranged
in a two-dimensional $2\times 4$ ladder with open boundary conditions, with spin up and down for each site.
A graphical illustration of the geometry is given in Figure~\ref{fig:geometry}.
Such an electron distribution is also referred to as \emph{half-filled} in the literature.
We furthermore restrict our model by considering the number of electrons with spin up or down to be fixed as $\nsit/2$,
respectively.
This leads to
$n=\genfrac{(}{)}{0pt}{1}{8}{4}^2=4900$ considered occupation states which create a discrete basis. Without restriction on the number of electrons and their spin, we would have $4^N= 65536$ states.
For the numerical implementation of the basis we consider $16$-bit integers
for which each bit describes a position which is occupied in case the bit is equal to $1$ or empty otherwise.
The set of occupation states can be ordered by the value of the integers which leads to a unique representation
of the Hubbard Hamiltonian~\eqref{eq.HubHam} by a matrix $H\in\C^{n\times n}$.
Such an implementation of the Hubbard Hamiltonian is also described in~\cite{Innerberger2020} and~\cite[Section 3]{Ja08}.
\item Secondly, we use a $4\times 3$ lattice ($\nsit=12$).
The state is again assumed as half-filled, that is, $12$ electrons populate the $12$ sites (6 with each spin).
This leads to $n=\genfrac{(}{)}{0pt}{1}{12}{6}^2=853776$ occupation states. The basis is also encoded by $16$-bit integers, analogously as in the $2\times 4$ case.~\cite{Innerberger2020}.
\end{enumerate}
\noindent In our first test setting, the $2\times 4$ ladder, we use $U=4$ and  a local potential $v_{11}=v_{88}=v_{44}=v_{55}=-1.75$ and $v_{22}=v_{33}=v_{66}=v_{77}=-2.25$ (see~Fig.\ref{fig:geometry}). Hopping is only allowed between nearest neighbors (with open boundary conditions). The absolute value of the nonzero $v_{i,j\neq i}$ elements is constant and  equal to $1$.
 This sets our energy units, whereas the unit of time $t$ is the inverse of the energy unit (which corresponds to setting $\hbar\equiv1$).

For this choice of $v_{ij}$ we obtain an Hermitian matrix $H(t)\in\C^{n\times n}$
with $60864$ nonzero entries for the $2\times 4$ ladder geometry.
The spectrum of $H(t)$ is independent of $t$ and the eigenvalues lie within the interval $(-21.04,5.23)$. This is because the electric field \eqref{eq:peierls} only modifies the phase of the nearest neighbor hopping elements according to
\begin{equation}\label{eq:pulse}
v_{ij}(t)=v_{ij}\underbrace{\exp\left( \ii  a (\cos(\omega(t-t_p))-b)\,\ee^{-\frac{(t-t_p)^2}{2\sigma_p^2}}\right)}_{=:f(t)},
\end{equation}
whereas the local terms $v_{ii}$ are constant in time.
Mathematically speaking, the described Hamiltonian $\Hfull(t)$
is \emph{isospectral}, since it satisfies \cite[Def. 8.3.8, \emph{Property A}]{stobul90}.
This also implies that in a Lanczos algorithm for the approximation of the
matrix exponential, the size of the subspace can be fixed for all $t$.
For our tests we choose $t_p=6$, $a=0.2$, $\sigma_p=2$, $\omega=3.5$, and $b=\cos(\omega t_p)$  for the $2\times4$ ladder. In supplementary comparisons given in \ref{further}
we also vary the parameters in the $2\times 4$ geometry and choose $\sigma_p=1,\ 4$ and $\omega=1.75,\ 7$
to corroborate our findings on the efficiency of the adaptive methods.

For the $4\times 3$ geometry, we use $U=8$, $v_{ii}=-4$ (the same for all sites), and the hopping elements $v_{ij}$ are again nonzero only for nearest neighbor sites $i$ and $j$ and given by~\eqref{eq:pulse} (which corresponds to a diagonal in-plane field). The parameters are as follows: $t_p=7.5$, $\sigma_p=2$, $a=0.8$, $\omega=11$, and $b=\cos(\omega t_p)$. In the $4\times3$ case we have $16687440$ nonzero entries of $H(t)$ and its eigenvalues lie within the interval $(-52.92, 4.92)$. For this set of parameters, impact ionization has been found~\cite{Kauch2020}.
%

Separating a constant diagonal contribution, which includes $v_{ii}$ and $U$, and splitting of the off-diagonal contribution into real and
imaginary part, this model leads to a Hamiltonian of the structure
\begin{equation}\label{pranavstructure}
\Hfull(t)=\Hdiag+\underbrace{{\rm Re}(f(t))}_{=:\,c(t)}\Hsymm
+\ii\underbrace{{\rm Im}(f(t))}_{=:\,s(t)}\Hanti
= \Hfull^\ast(t),
\end{equation}
with real matrices $ \Hdiag $, $ \Hsymm $, and $ \Hanti $,
where $ \Hdiag $ is diagonal, $ \Hsymm $ is symmetric, and
$ \Hanti $ is skew-symmetric.
Note that the time-dependence is only present in the scalar function
$f(t)$, while the matrices are constant.
This structure will be exploited in a new time integrator using
operator splitting, which we introduce in Section~\ref{subsec:MS}, see also
\ref{sec:adapt}.

\section{Numerical approaches} \label{sec:num-approach}

We consider Magnus-type one-step methods for
the approximation of~\eqref{eq0} on a time grid $(t_0,t_1,\dots,t_n,\dots)$,
\begin{equation*}
\psi_{n+1} = \nS(\t_n;t_n)\,\psi_n \approx \psi(t_{n+1})
           = \nE(\t_n;t_n)\,\psi(t_n),\quad \t_n=t_{n+1}-t_n\,,
\quad n=0,1,2,\ldots\,,
\end{equation*}
where $\nE(\tau_n;t_n)$ denotes the exact and $\nS(\tau_n;t_n)$ the approximated unitary operator for the time propagation from $t_n$ to $t_n+\tau_n$.
For the description of the schemes,
in the following we use a simplified notation for a single step
starting from $ t=t_0 $ with stepsize $ \t $,
\begin{equation}
\label{one-step-0}
\psi_{1} = \nS(\t;t_0)\,\psi_0 \approx \psi(t_0+\t)\,.
\end{equation}
In order to avoid unnecessary overloading of notation,
we suppress the dependence on $ t_0 $ of `internal' objects
involved in the definition of the integrators.
Only the dependence on the stepsize $ \t $ is indicated;
see for instance~\eqref{CFM} below.

In the following, we will put our emphasis on approximations of order four,
as this is usually sufficient for the accuracy required, and
the splitting-based integrator described in Section~\ref{subsec:MS} has this order.
All the exponential-based time integrators considered in this study are
symmetric (time-reversible), which is a desirable property when
the reversible flow of a Schr{\"o}dinger-type equation is approximated.

\subsection{Commutator-free Magnus-type (CFM) integrators}\label{subsec:CFM}

A successful and much used class of integration methods is comprised of
higher-order commutator-free Magnus-type integrators~\cite{alvfeh11,blamoa05}.
These approximate the exact flow in terms of
products of exponentials of linear combinations of the system matrix evaluated at different times,
avoiding evaluation and storage of commutators.

A high-order CFM scheme starting at $t=t_0$ is thus defined by~\eqref{one-step-0}, with the ansatz
\cite{alvfeh11,blamoa05}
\begin{equation} \label{CFM}
\begin{aligned}
\nS(\t;t_0) &=
\nS_J(\t) \cdots \nS_1(\t)
= \ee^{\Omega_J (\t)} \,\cdots\, \ee^{\Omega_1(\t)}\,, \\[1mm]
\Omega_j(\t) &= \t B_j(\t),~~ j=1, \ldots, J, \\
B_{j}(\t) &= \sum_{k=1}^K a_{jk}\, A_{k}(\t),\qquad A_{k}(\t) = A(t_0+c_k \t)\,,
\end{aligned}
\end{equation}
where the coefficients $a_{jk}$, $c_k$ are determined from the \emph{order conditions}
(a system of polynomial equations in the coefficients) such that the method realizes a certain
convergence order $ p $, see for example \cite{SergioFernandoMPaper2} and references therein.
Algorithms to efficiently generate the order conditions are described for instance in \cite{moskau19}.
The solution of this system of equations is generally not unique, and numerical optimization
techniques are employed to compute solutions that are optimal in some sense, for instance
minimizing in some sense the \emph{leading local error term} of the ensuing integrator.

\paragraph{Examples of symmetric CFM integrators}
\begin{enumerate}[(i)]
\item
The second-order scheme ($ p=2 $) given by
\begin{equation*}
J = 1, \quad K = 1, \qquad c = \tfrac{1}{2}, \quad a = 1\,,
\end{equation*}
is a simple instance of a Magnus-type integrator, commonly
denoted as \emph{exponential midpoint rule}. It yields
\begin{equation}\label{CF2}
\nS(\t;t_0) = \ee^{\t A(t_0+\frac{\t}{2})}\,.
\end{equation}
Note that this represents both a classical Magnus integrator and a commutator-free
Magnus-type method. In the numerical tests, this method is referred to as \texttt{CF2}.
\item
A fourth-order commutator-free integrator ($ p=4 $)
based on two Gaussian nodes and
comprising two matrix exponentials is defined by $J=K=2$ and
\begin{equation}
\label{CF4}
c = \begin{pmatrix}
\tfrac{1}{2} - \tfrac{\sqrt{3}}{6} \\
\tfrac{1}{2} + \tfrac{\sqrt{3}}{6}
\end{pmatrix},
\quad
a = \begin{pmatrix}
    \frac{1}{4}+\frac{\sqrt{3}}{6} & \frac{1}{4}-\frac{\sqrt{3}}{6} \nl{1}
    \frac{1}{4}-\frac{\sqrt{3}}{6} & \frac{1}{4}+\frac{\sqrt{3}}{6}
    \end{pmatrix}.
\end{equation}
In the numerical tests, this method is referred to as \texttt{CF4}.
\item An optimized fourth-order scheme ($ p=4 $) from~\cite{alvfeh11}
satisfies $J=K=3$ and 
\begin{equation} \label{CF-4o}
c =
\begin{pmatrix}
\frac{1}{2}-\frac{\sqrt{15}}{10} \nl{2}
\frac{1}{2}                      \nl{1}
\frac{1}{2}+\frac{\sqrt{15}}{10}
\end{pmatrix},
\quad
a =
\begin{pmatrix}
  \frac{37}{240} + \frac{10}{87}\frac{\sqrt{15}}{3} &
 -\frac{1}{30}                                     &
  \frac{37}{240} - \frac{10}{87}\frac{\sqrt{15}}{3}     \nl{2}
 -\frac{11}{360}                                   &
 ~\;\frac{23}{45}                                    &
 -\frac{11}{360}                                       \nl{1}
  \frac{37}{240} - \frac{10}{87}\frac{\sqrt{15}}{3} &
 -\frac{1}{30}                                     &
  \frac{37}{240} + \frac{10}{87}\frac{\sqrt{15}}{3}
\end{pmatrix}.
\end{equation}
In the numerical tests, this method is referred to as \texttt{CF4o}.
\item Based on the approach for constructing new optimized commutator-free
Magnus-type integrators described in \cite{moskau19,hofi19}, we have constructed the following
fourth-order numerical integrator with $J=K=3,$ but smaller leading error term:
\begin{eqnarray*} \label{CF-4oH}
&&c =
\begin{pmatrix}
1/2-\sqrt{15}/10  \nl{2}
1/2 \nl{1}
2+\sqrt{15}/10
\end{pmatrix},
\\[2mm]
&&a_{11}=0.302146842308616954258187683416\\
&&a_{12}=-0.030742768872036394116279742324\\
&&a_{13}=0.004851603407498684079562131338\\
&&a_{21}=-0.029220667938337860559972036973\\
&&a_{22}=0.505929982188517232677003929089\\
&&a_{23}=-0.029220667938337860559972036973\\
&&a_{31}=0.004851603407498684079562131337\\
&&a_{32}=-0.030742768872036394116279742324\\
&&a_{33}=0.302146842308616954258187683417
\end{eqnarray*}

In the numerical tests, this method is referred to as \texttt{CF4oH}.
Note that the coefficients have been determined numerically and are given to within double precision,
likewise for the next schemes.
\item We have also constructed a new optimized method of order six with $K=3,\ J=4$:
\begin{eqnarray*} \label{CF6n}
&& c =
\begin{pmatrix}
1/2-\sqrt{15}/10 \nl{2}
1/2 \nl{1}
1/2+\sqrt{15}/10
\end{pmatrix},
\\[2mm]
&&a_{11}=0.79124225942889763\\
&&a_{12}=-0.080400755305553218\\
&&a_{13}=0.01.2765293626634554\\
&&a_{21}=-0.48931475164583259\\
&&a_{22}=0.05.4170980027798808\\
&&a_{23}=-0.012069823881924156\\
&&a_{31}=-0.029025638294289255\\
&&a_{32}=0.50138457552775674\\
&&a_{33}=-0.02.5145341733509552\\
&&a_{41}=0.0048759082890019896\\
&&a_{42}=-0.030710355805557892\\
&&a_{43}=0.30222764976657693
\end{eqnarray*}

In the numerical tests, this method is referred to as \texttt{CF6n}.
\item A new optimized method of order seven with $K=4,\ J=6$ is given by the coefficients
\begin{eqnarray*} \label{CF7}
&& c =
\begin{pmatrix}
-\sqrt{\tfrac{1}{140}(2\sqrt{30}+15)}+1/2\nl{3}
-\sqrt{\tfrac{1}{140}(-2\sqrt{30}+15)}+1/2\nl{2}
\sqrt{\tfrac{1}{140}(-2\sqrt{30}+15)}+1/2\nl{1}
\sqrt{\tfrac{1}{140}(2\sqrt{30}+15)}+1/2
\end{pmatrix},
\\[2mm]
&&a_{11}=0.205862188450411892209\\
&&a_{12}=0.169508382914682544509\\
&&a_{13}=-0.102088008415028059851\\
&&a_{14}=0.0304554010755044437431\\
&&a_{21}=-0.0574532495795307023280\\
&&a_{22}=0.234286861311879288330\\
&&a_{23}=0.332946059487076984706\\
&&a_{24}=-0.0703703697036401378340\\
&&a_{31}=-0.008.93040281749440468751\\
&&a_{32}=0.02.71488489365780259156\\
&&a_{33}=-0.02.95144169823456538040\\
&&a_{34}=-0.151311830884601959206\\
&&a_{41}=0.552299810755465569835\\
&&a_{42}=-3.64425287556240176808\\
&&a_{43}=2.53660580449381888484\\
&&a_{44}=-0.661436528542997675116\\
&&a_{51}=-0.538241659087501080427\\
&&a_{52}=3.60578285850975236760\\
&&a_{53}=-2.50685041783117850901\\
&&a_{54}=0.651947409253201845106\\
&&a_{61}=0.02.03907348473756540850\\
&&a_{62}=-0.0664014986792173869631\\
&&a_{63}=0.0949735566789294244299\\
&&a_{64}=0.374643341371260411994
\end{eqnarray*}
In the numerical tests, this method is referred to as \texttt{CF7}.
\end{enumerate}

\paragraph{Local error estimation}

As a basis for adaptive time-stepping, defect-based error estimators for
CFM methods and for classical Magnus integrators have been introduced in \cite{auzingeretal18b}.
The (classical) defect is defined by
\begin{equation}\label{defect1}
\nD(\tau)=\nS'(\tau;t_0)-A(t_0+\tau)\nS(\tau;t_0)
\end{equation}
and satisfies
$$\nD(0)=\nD'(0)=\cdots=\nD^{(p-1)}(0)=0,$$
if the method coefficients $a,\ c$ satisfy the order conditions for an order $p$ method.

The local error $\nL(\tau)\psi_0:=(\nS(\tau;t_0)-\nE(\tau;t_0))\psi_0$ can be expressed in terms
of the defect via the \emph{variation-of-constant formula},
$$\nL(\tau)\psi_0=\int_0^\tau\Pi(\tau,\sigma)\nD(\sigma)\,\dd\sigma=\Order(\tau^{p+1}),\quad
\Pi(\tau,\sigma)=\nE(\tau-\sigma;t_0+\sigma).$$
For the practical evaluation of the defect, the derivative of matrix exponentials of the form
$$\dt \ee^{\tau B(\tau)} = \Gamma(\tau)\,\ee^{\tau B(\tau)}$$
is required.
The function $\Gamma$ is given as an infinite series or alternatively as an integral
expression, which are approximated by truncation or numerical Hermite quadrature, respectively,
to yield a computable quantity $\tilde\Gamma$ and an approximate defect $\tilde\nD$.
The resulting computable error estimator is denoted by $\tilde\nP$ in both cases.
The \emph{asymptotical correctness} of the error estimators was established in \cite{auzingeretal18b}.
We recapitulate the result for the exponential midpoint rule (\ref{CF2}):

\medskip

\noindent\textbf{Proposition:} \textit{Consider the exponential midpoint rule {\rm (\ref{CF2})}.\\
If $A\in C^3$, then the local error $\nL$ satisfies
$$
\| \nL(\t) \|_2 \leq
\tfrac{1}{12} \t^3 \big\| [A(t_0),A'(t_0)]
                         - \tfrac{1}{2} A''(t_0) \|_2
+ \Order(\t^4)\,.
$$
If $A\in C^4$, then the deviation
$ \tilde\nP(\t)-\nL(\t) =
 \tfrac{1}{3}\t\tilde\nD(\t)-\nL(\t) $
of the local error estimate satisfies
\begin{eqnarray}
\| \tilde\nP(\t)-\nL(\t) \|_2
&\leq& \tau^4 \big\|c\,[A(t_0),[A(t_0),A'(t_0)]] - \tfrac{1}{48}[A(t_0),A''(t_0)]\nonumber \\
&& \qquad +\tfrac{1}{144}A'''(t_0)\big\|_2 + \Order(\t^5)\,, \label{deviation1}
\end{eqnarray}
where {$ c=\tfrac{1}{72} $} for the approximate defect $ \tilde\nD(\t) $, Taylor version,
and {$ c=0 $} for the approximate defect $ \tilde\nD(\t) $, Hermite version.}

In the numerical experiments reported in this paper, Hermite quadrature has been used throughout.

\subsection{Classical Magnus integrators}
A different, indeed the more classical,
approach to the approximation of~\eqref{eq0}
is directly based on the Magnus expansion~\cite{magnus54}:
The solution to a time-dependent system~\eqref{eq0} can be represented by
\begin{subequations} \label{magnus1}
\begin{equation} \label{magnus1a}
\psi(t_0+\t) = \nE(\t;t_0)\psi_0 = \ee^{\bm\Omega(\t)} \psi_0\,,
\end{equation}
where the matrix $\bm\Omega(\t)$ satisfies
\begin{equation}\label{magnus1b}
\bm\Omega'(\t)
= \sum_{k \geq 0}\frac{ {B_k}}{k!}\,\mathrm{ad}_{\bm\Omega(\t)}^k(A(t_0+\t))\,,
\quad \bm\Omega(0) = 0\,,
\end{equation}
\end{subequations}
with the Bernoulli numbers $B_k$ and $\mathrm{ad}_{\bm\Omega}^k(A)$ denoting
the $k$-th iterated commutator of the matrices ${\bm\Omega}$ and $A$, see \cite{haireretal02b}.

Classical Magnus integrators rely on appropriate truncation of
the Magnus expansion~\eqref{magnus1b} and
suitable approximation $ \Omega(\t) $ to the arising multi-dimensional
integral representation for $ \bm\Omega(\t) $ by numerical quadrature,
and defining $ \psi_1 $ by~\eqref{one-step-0} with
\begin{equation} \label{magnus-classic}
\nS(\t;t_0) = \ee^{\Omega(\t)}\approx\ee^{\bm\Omega(\t)}.
\end{equation}
A detailed exposition on this approach is given for example
in~\cite{blanesetal08b} and in~\cite{haireretal02b},
see also \cite{iserlesetal00}.

This type of integrator is, in general, considered as computationally expensive
due to the requirement to compute and store commutators of large matrices.
For problems of a particular structure, however, as in the semiclassical regime,
where  the small semiclassical parameter may render some of the appearing
commutators negligibly small,
or when commutators turn out to be of higher order $O(\t^k)$ than
$O(1)$ as expected generically,
this approach may excel over the commutator-free methods,
see~\cite{alvfeh11,SergioFernandoMPaper2,baderetal17}. In our specially designed numerical
integrator in Section~\ref{subsec:MS} below, such a feature is actually exploited, see~\cite{karo18a}.

\paragraph{Examples of classical symmetric Magnus integrators}
\begin{enumerate}[(i)]
\item The exponential midpoint scheme~\eqref{CF2} (order $ p=2 $)
is also a classical Magnus integrator, with $ K=1 $ and
\begin{equation}\label{M2}
c = \tfrac{1}{2}\,, \quad \Omega(\t) = \t A_1(\t)\,.
\end{equation}
\item
A commonly used fourth-order Magnus integrator ($ p=4 $)
is based on two Gaussian nodes, with $ K=2 $ and
\begin{equation}\label{M4}
c =
\begin{pmatrix}
\tfrac{1}{2} - \tfrac{\sqrt{3}}{6}\\
\tfrac{1}{2} + \tfrac{\sqrt{3}}{6}
\end{pmatrix},
\quad
\Omega(\t) = \tfrac{1}{2} \t \big(A_1(\t) + A_2(\t)\big)
              - \tfrac{\sqrt{3}}{12}\,\t^2 \big[A_1(\t), A_2(\t)\big]\,.
\end{equation}
In the numerical experiments, this is denoted as \texttt{Magnus4}.
\end{enumerate}

\paragraph{Local error estimation}
Classical Magnus integrators are of the form~\eqref{magnus-classic},
where again $ \Omega(\t) = \t B(\t) $. Thus, the defect
can be approximated analogously as in Section~\ref{subsec:CFM}.
Convergence of the integrators and asymptotical correctness of the error estimator has been proven
in \cite{auzingeretal18b}.

\subsection{Symmetrized defect}\label{sec:symmdef}

When we consider \emph{self-adjoint} (or \emph{symmetric}\/) schemes which are
characterized by the identity
\begin{equation} \label{nonaut-symmscheme}
\nS(-\tau;t_0+\tau)\,\nS(\tau;t_0) = \Id,
\end{equation}
a higher asymptotical quality of the error estimator can be obtained at moderate additional expense
by introducing a \emph{symmetrized version of the defect},
which was introduced and analyzed in \cite{auzingeretal18a,auzingeretal18c}. We define
\begin{align}
&\nDs(\tau) =
\nS'(\tau;t_0) -
\th\big( A(t_0+\tau) \nS(\tau;t_0)
          + \pdtwo\nS(\tau;t_0)
           + \nS(\tau;t_0)A(t_0)
   \big), \label{dst}
\end{align}
where $\pdtwo$ stands for differentiation with respect to the second argument.
A local error representation based on a symmetrized variation-of-constant formula
in conjunction with numerical quadrature again yields an asymptotically correct
error estimator. Since the exponential-based integrators we used in this study are symmetric,
the symmetric version of the error estimator has been used throughout.

\paragraph{Example: Exponential midpoint rule \cite{auzingeretal18c}}\label{subsec:emr}
Let
\begin{equation*}
\nR(\tau;t_0)(\,\cdot\,) =
\tfrac{\dd}{\dd\Omega} \ee^{\Omega}\big|_{\Omega=\tau A(t_0+\frac{\tau}{2})}(\,\cdot\,),
\end{equation*}
where $ \tfrac{\dd}{\dd\Omega} \ee^\Omega $ denotes the Fr{\'e}chet
derivative of the matrix exponential, see~\eqref{mimi} below. Then,
for the exponential midpoint rule~\eqref{CF2}
\begin{align*}
\nS'(\tau;t_0)
&= \nR(\tau;t_0) \big( A(t_0+\tfrac{\tau}{2}) + \th\tau A'(t_0+\tfrac{\tau}{2}) \big), \\
\pdtwo \nS(\tau;t_0)
&= \nR(\tau;t_0) \big( \tau A'(t_0+\tfrac{\tau}{2}) \big).
\end{align*}
This gives the following defect representations:
\begin{itemize}
\item
Classical defect~\eqref{defect1}:
\begin{equation} \label{dcA(t)}
\nDc(\tau) =
\nR(\tau;t_0) \big( A(t_0+\tfrac{\tau}{2})
                    + \th\tau A'(t_0+\tfrac{\tau}{2}) \big)
- A(t_0+\tau) \nS(\tau;t_0).
\end{equation}
\item
Symmetrized defect~\eqref{dst}:
\begin{align}
&\nDs(\tau) =
\nR(\tau;t_0) \big( A(t_0+\tfrac{\tau}{2}) \big)
   - \th \big( A(t_0+\tau) \nS(\tau;t_0) + \nS(\tau;t_0) A(t_0) \big). \label{dsA(t)}
\end{align}
\end{itemize}
Here, the explicit representation
\begin{align} \label{mimi}
\nR(\tau;t_0) \big( V \big)
&= \int_0^1 \ee^{\sigma\tau A(t_0+\frac{\tau}{2})}
            V
            \ee^{-\sigma\tau A(t_0+\frac{\tau}{2})}\,\dd\sigma
   \cdot \nS(\tau;t_0)
\end{align}
follows from~\cite[(10.15)]{higham08}.
For evaluating~\eqref{dcA(t)}, a sufficiently accurate quadrature
approximation for the integral according to~\eqref{mimi} is required.
This involves evaluation of $ A' $ and the commutator $ [A,A'] $,
see~\cite{auzingeretal18c}.
In contrast, the relevant term from~\eqref{dsA(t)} simplifies to
\begin{align*}
\nR(\tau;t_0) \big( A(t_0+\tfrac{\tau}{2}) \big)
&= A(t_0+\tfrac{\tau}{2}) \nS(\tau;t_0)
 = \nS(\tau;t_0) A(t_0+\tfrac{\tau}{2}),
\end{align*}
whence the symmetrized defect~\eqref{dsA(t)} can be evaluated exactly,
\begin{equation} \label{dsA(t)-1}
\begin{aligned}
\nDs(\tau)
&= \nS(\tau;t_0) \big( A(t_0+\tfrac{\tau}{2}) - \th A(t_0) \big)
  - \th A(t_0+\tau) \nS(\tau;t_0).
\end{aligned}
\end{equation}
This involves an additional application of $ \nS(\tau;t_0) $,
but it does not require evaluation of the derivative $ A' $
or of a commutator expression.
We also note that the applications of
$\nS$ from left and right can be evaluated in parallel.

The deviation of the symmetrized error estimator used in conjunction
with the exponential midpoint rule as compared to
the exact local error satisfies an estimate of the form
$$\| \tilde\nP(\t)-\nL(\t) \|_2 =\Order(\t^5),$$
see~\cite{auzingeretal18c}.

\subsection{Fourth order Magnus-Strang splitting} \label{subsec:MS}
\paragraph{Magnus expansion} To obtain a new efficient fourth-order method
for problems of the structure (\ref{pranavstructure}),
we proceed similarly as in \cite{karo18a}.
We start using only the first two terms of the Magnus expansion \eqref{magnus1}
for the solution of~\eqref{eq0},
$$
\psi(t_0+\tau) = {\rm e}^{\mathbf{\Omega}(\tau)} \psi_0 \approx {\rm e}^{\Omega_2(\tau)} \psi_0 + \mathcal{O}(\tau^5),
$$
with
\begin{equation} \label{magnus-trunc}
\begin{split}
\Omega_2(\tau)
&= \int_{t_0}^{t_0+\tau} -\ii\,\Hfull(\zeta)\,{\rm d}\zeta
+ \frac12\int_{t_0}^{t_0+\tau}\int_{t_0}^{t_0+\zeta}[-\ii\,\Hfull(\zeta),-\ii\,\Hfull(\xi)]\,{\rm d}\xi\,{\rm d}\zeta \\
&= \underbrace{\Big(-\ii \int_{t_0}^{t_0+\tau} \Hfull(\zeta)\,{\rm d}\zeta\Big)}_{\nA(\tau;t_0)=\mathcal{O}(\tau)}
+ \underbrace{\Big(-\frac12\int_{t_0}^{t_0+\tau}\int_{t_0}^{t_0+\zeta}[\Hfull(\zeta),\Hfull(\xi)]\,{\rm d}\xi\,{\rm d}\zeta\Big)}_{\nB(\tau;t_0)=\mathcal{O}(\tau^3)},
\end{split}
\end{equation}
where $ \nA(\tau;t_0) $ and $ \nB(\tau;t_0) $ are again skew-Hermitian.
Here, because the time-dependence in (\ref{pranavstructure}) is given by two scalars  $s(t)$ and $c(t)$, $ \nA(\tau;t_0) $ is simple to evaluate,
$$
\int_{t_0}^{t_0+\tau}\Hfull(\xi)\,{\rm d}\xi=t\Hdiag+\Hsymm\int_{t_0}^{t_0+\tau} c(\xi)\,{\rm d}\xi+\ii \Hanti\int_{t_0}^{t_0+\tau} s(\xi)\,{\rm d}\xi,
$$
while evaluation of $ \nB(\t;t_0) $ at first sight seems more challenging.
It involves evaluation of the commutator
\begin{align}
[\Hfull(\zeta),\Hfull(\xi)]
&=c(\xi)[\Hdiag,\Hsymm]+\ii\,s(\xi)[\Hdiag,\Hanti]\nonumber\\
& \quad {} + c(\zeta)[\Hsymm,\Hdiag]+\ii\,s(\zeta)[\Hanti,\Hdiag]\nonumber\\
& \quad {} + \ii\,c(\zeta)s(\xi)[\Hsymm,\Hanti]+\ii\,c(\xi)s(\zeta)[\Hanti,\Hsymm]\nonumber\\
&=\left(c(\zeta)-c(\xi)\right)[\Hsymm,\Hdiag]+\ii \left(s(\zeta)-s(\xi)\right)[\Hanti,\Hdiag]\nonumber\\
& \quad {} +\ii \left(c(\zeta)s(\xi)-c(\xi)s(\zeta)\right)[\Hsymm,\Hanti].  \label{Hfull2}
\end{align}
Since $\Hdiag$, $\Hsymm$, and $\Hanti$ do not depend on time,
the computation of $\int_{t_0}^{t_0+\tau}\int_{t_0}^{t_0+\zeta}[\Hfull(\zeta),\Hfull(\xi)]\,{\rm d}\xi\,{\rm d}\zeta$
boils down to the evaluation of the following integrals, which can be simplified as indicated:
\begin{equation} \label{integrals-simplified}
\begin{split}
&\int_{t_0}^{\tau}\int_{t_0}^{t_0+\zeta}\big(c(\zeta)-c(\xi)\big)\,{\rm d}\xi\,{\rm d}\zeta
= 2\int_{t_0}^{t_0+\tau}c(\zeta)\big(\zeta-\tfrac{\tau}{2}\big)\,{\rm d}\zeta, \\
&\int_{t_0}^{t_0+\tau}\int_{t_0}^{t_0+\zeta} \big(s(\zeta)-s(\xi)\big)\,{\rm d}\xi\,{\rm d}\zeta
= 2 \int_{t_0}^{t_0+\tau}s(\zeta)(\zeta-\tfrac{\tau}{2})\,{\rm d}\zeta, \\
&\int_{t_0}^{t_0+\tau}\int_{t_0}^{t_0+\zeta} \big(c(\zeta)s(\xi)-c(\xi)s(\zeta)\big)\,{\rm d}\xi\,{\rm d}\zeta \\
& \quad = \int_{t_0}^{t_0+\tau} c(\zeta)\,{\rm d}\zeta \int_{t_0}^{t_0+\tau} s(\xi)\,{\rm d}\xi-2\int_{t_0}^{t_0+\tau}c(\zeta)\int_{t_0+\zeta}^{t_0+\tau} s(\xi)\,{\rm d}\xi\,{\rm d}\zeta.
\end{split}
\end{equation}

\paragraph{The splitting}
To obtain a practical fourth order method
we proceed in a similar way as indicated in~\cite{karo18a}.
Applying Strang splitting to~\eqref{magnus-trunc} yields
$$
{\rm e}^{\frac12 \nA(\tau;t_0)} {\rm e}^{\nB(\tau;t_0)} {\rm e}^{\frac12 \nA(\tau;0)}
= {\rm e}^{\mathbf{\Omega}(\tau;t_0)} + \mathcal{O}(\tau^5),
$$
or alternatively
$$
{\rm e}^{\frac12 \nB(\tau;t_0)} {\rm e}^{\nA(\tau;t_0)} {\rm e}^{\frac12 \nB(\t;t_0)}
= {\rm e}^{\mathbf{\Omega}(\tau)} + \mathcal{O}(\tau^5).
$$
Which of the two splitting variants is favorable depends on the properties of the operators $\nA$ and $\nB$.
For the present problem the splitting can be realized in an efficient way, expressing $ \nA $ and $ \nB $
in the form
\begin{align*}
\nA(\tau;t_0) &= -\ii\,\big( \tau \Hdiag+\tilde{c}_1(\tau,t_0)\Hsymm+\ii\,\tilde{s}_1(\tau,t_0)\Hanti \big) \\
     &= -\ii\,\tau \Hdiag-\ii\,\tilde{c}_1(\tau,t_0)\Hsymm+\tilde{s}_1(\tau,t_0)\Hanti,
\end{align*}
and (see~\eqref{Hfull2})
\begin{equation*}
\nB(\tau;t_0) = -\tilde{c}_2(\tau,t_0)[\Hsymm,\Hdiag]-\ii\,\tilde{s}_2(\tau,t_0)[\Hanti,\Hdiag]-\ii\,\tilde{r}(\tau,t_0)[\Hsymm,\Hanti],
\end{equation*}
where (see~\eqref{integrals-simplified})
\begin{align*}
\tilde{c}_1(\tau;t_0)&=\int_{t_0}^{t_0+\tau}c(\zeta)\,{\rm d}\zeta= \mathcal{O}(\tau),\quad
\tilde{s}_1(\tau;t_0)=\int_{t_0}^{t_0+\tau}s(\zeta)\,{\rm d}\zeta=\mathcal{O}(\tau), \\
\tilde{c}_2(\tau;t_0)&=\int_{t_0}^{t_0+\tau}c(\zeta)(\zeta-\tfrac{\tau}{2})\,{\rm d}\zeta=\mathcal{O}(\tau^3),\\
\tilde{s}_2(\tau;t_0)&=\int_{t_0}^{t_0+\tau}s(\zeta)(\zeta-\tfrac{t}{2})\,{\rm d}\zeta=\mathcal{O}(\tau^3),\\
\tilde{r}(\tau;t_0)&=\tfrac{1}{2}\int_{t_0}^{t_0+\tau} c(\zeta)\,{\rm d}\zeta \int_{t_0}^{t_0+\tau} s(\xi)\,{\rm d}\xi-\int_{t_0}^{t_0+\tau}c(\zeta)\int_{t_0+\zeta}^{t_0+\tau} s(\xi)\,{\rm d}\xi\,{\rm d}\zeta=\mathcal{O}(\tau^3),
\end{align*}
and apply Strang splitting as before. Implementation details for this integrator are explained in \ref{sec:adapt}.
In the numerical experiments, this method is denoted as \texttt{MagnusStrang4}.

\begin{remark}
Methods based on a Strang splitting of the Magnus expansion were found to be very effective
when the integrand is highly oscillatory \cite{karo18a,pranav19}.
Note that due to the special scalings of $\nA$ and  $\nB$, this splitting is actually of order four \cite{karo18a}.

\end{remark}

\subsection{Adaptive Lanczos method}\label{sec:lanczos}

In each step of any of the introduced  Magnus-type methods, the action of a matrix exponential
\begin{equation}\label{expmat}
E(t)v=\ee^{-\ii t \Omega}v,\quad \text{$\Omega$\ Hermitian},
\end{equation}
has to be approximated. The standard Krylov approximation is
\begin{equation}\label{krylov}
S_m(t)v  = V_m\,\ee^{-\ii t T_m}\,V_m^\ast\,v
         = V_m\,\ee^{-\ii t T_m}e_1,
\end{equation}
with $T_m=(\tau_{i,j})$ tridiagonal and $V_m$ an orthonormal basis of the Krylov space
$\mathcal{K}_m(\Omega,v) = \mathrm{span}\{v,\Omega v,\ldots,\Omega^{m-1}v\} \subseteq \C^n$.
For Hermitian or skew-Hermitian matrices $\Omega$, this can be realized cheaply by the
Lanczos method \cite{molloa03}.
In \cite{jaweckietal20}, a time-stepping strategy was introduced which is based
on the defect of the approximation. Due to the success of this strategy
documented ibidem, we use it invariantly throughout our numerical experiments.
The asymptotically correct error estimator is described in the following:

First, we define the \emph{defect operator}
$$ D_m(t) = -\ii\,\Omega\,S_m(t) - S_m'(t) \in \C^{n \times n}.$$
Then the local error operator $L_m(t)=E(t)-S_m(t)$ enjoys the representation
$$ L_m(t)v = \int_0^t E(t-s)\,D_m(s)v\,\dd s.$$
This defect-based integral representation yields a computable,
asymptotically correct local error bound ($\gamma_m
          = \prod_{j=1}^{m-1} (T_m)_{j+1,j})$, satisfying (see \cite{jaweckietal20}).
\begin{eqnarray*}
&&\| L_m(t)v \|_2 \leq \tau_{m+1,m}\gamma_m\,\frac{t^{m}}{m!},\\[3mm]
&&\| L_m(t)v \|_2 = \tau_{m+1,m}\gamma_m\,\frac{t^m}{m!} + {\mathcal{O}}(t^{m+1}).
\end{eqnarray*}

\section{Numerical experiments} \label{sec:numexp}

In this section, we give the results of our numerical experiments to assess
the accuracy, reliability and efficiency of the numerical time
integrators investigated in this work. In all the experiments,
the tolerance for the adaptive Lanczos method for the matrix
exponential was set to $10^{-12}$. This limits the accuracy
that can be obtained by the time integrators, which manifests
itself as an apparent order reduction for the smallest time-steps in some cases.
Supplementary tests given in \ref{further} show similar results.

\begin{figure}[ht!]
\begin{center}
\includegraphics[width=5.8cm]{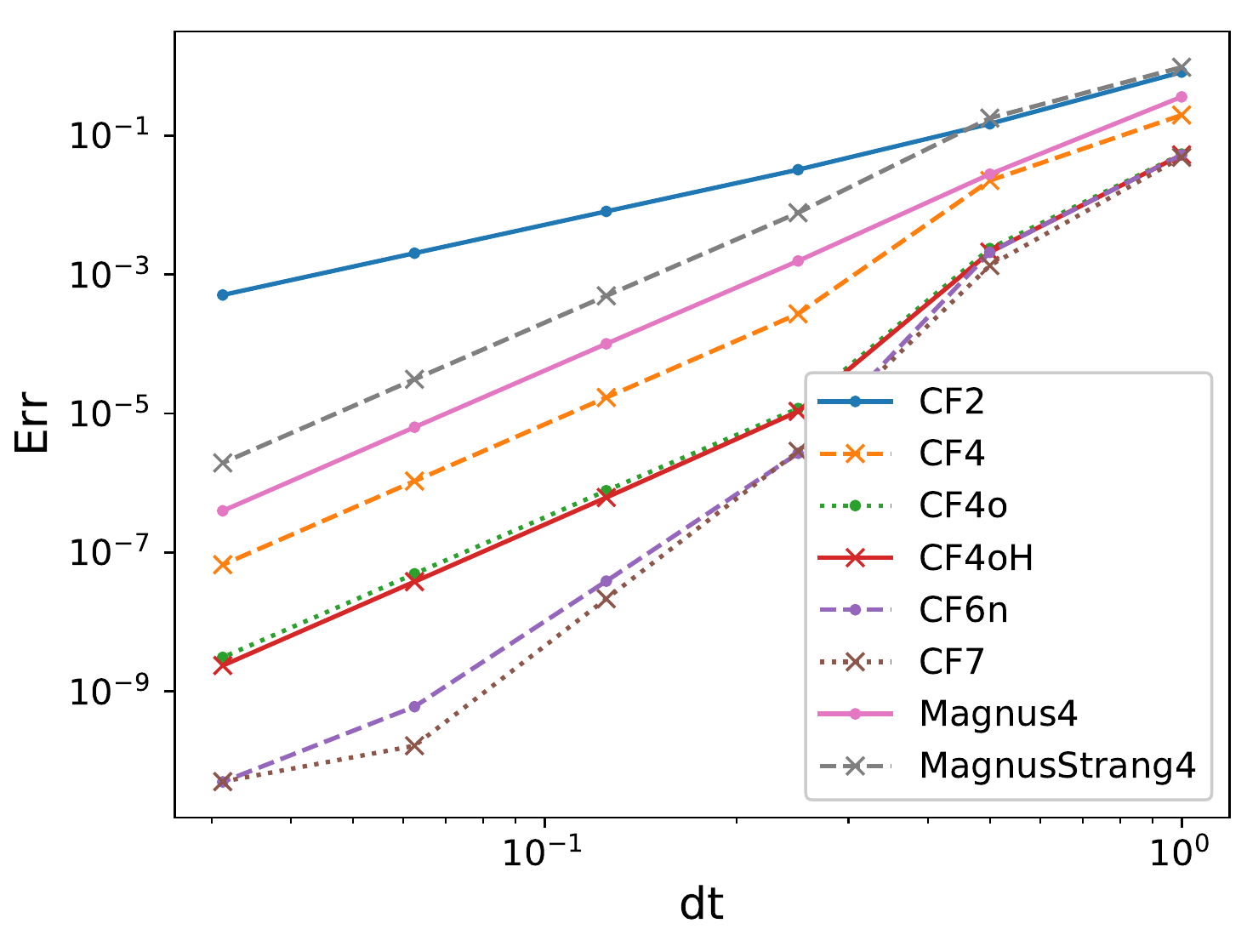} \quad \includegraphics[width=5.8cm]{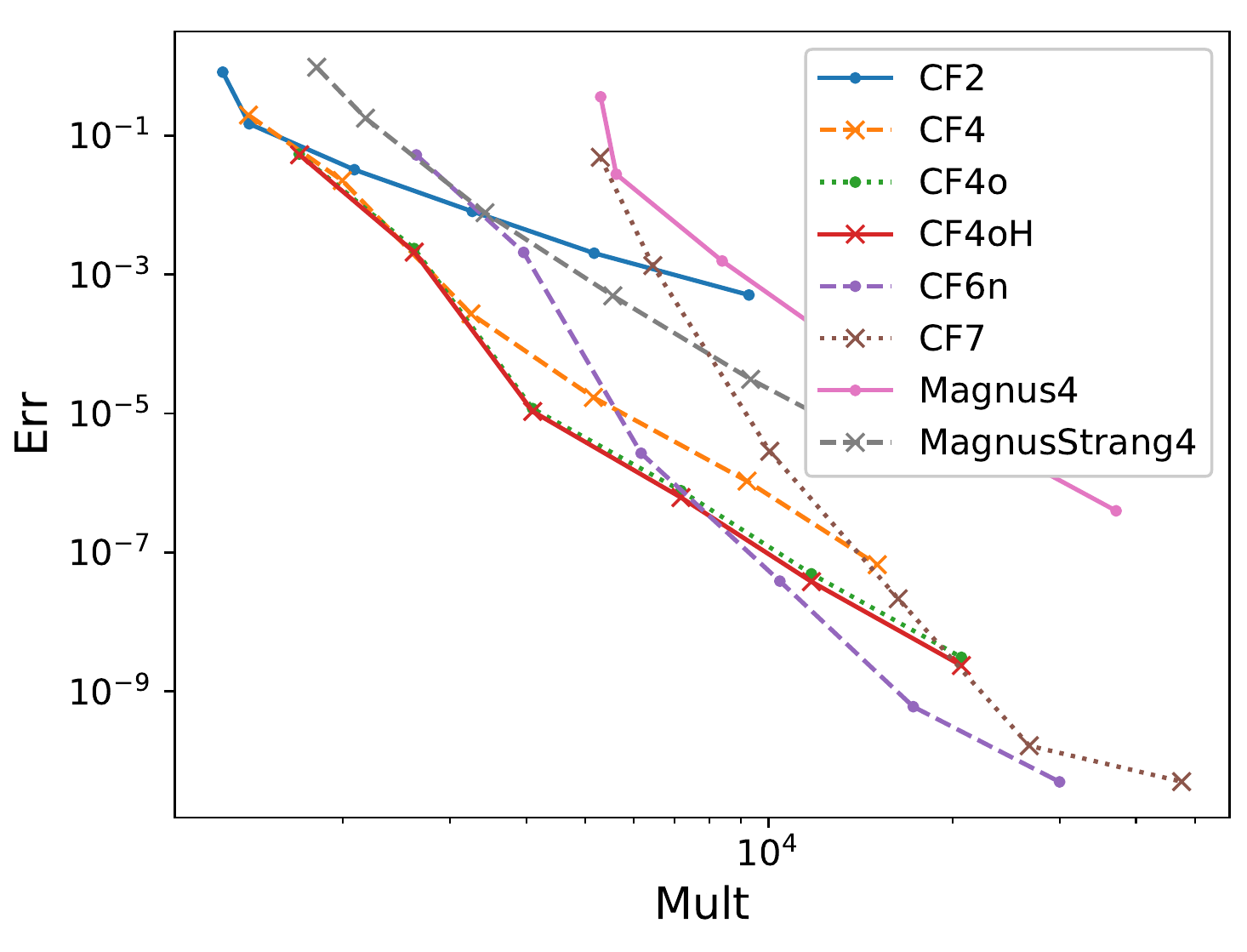}
\caption{ $2\times4$ geometry, equidistant time-steps. Error as a function of the step-size (left) and as
a function of matrix-vector multiplications (right).\label{equi8x1}}
\end{center}
\end{figure}

\begin{figure}[ht!]
\begin{center}
\includegraphics[width=5.8cm]{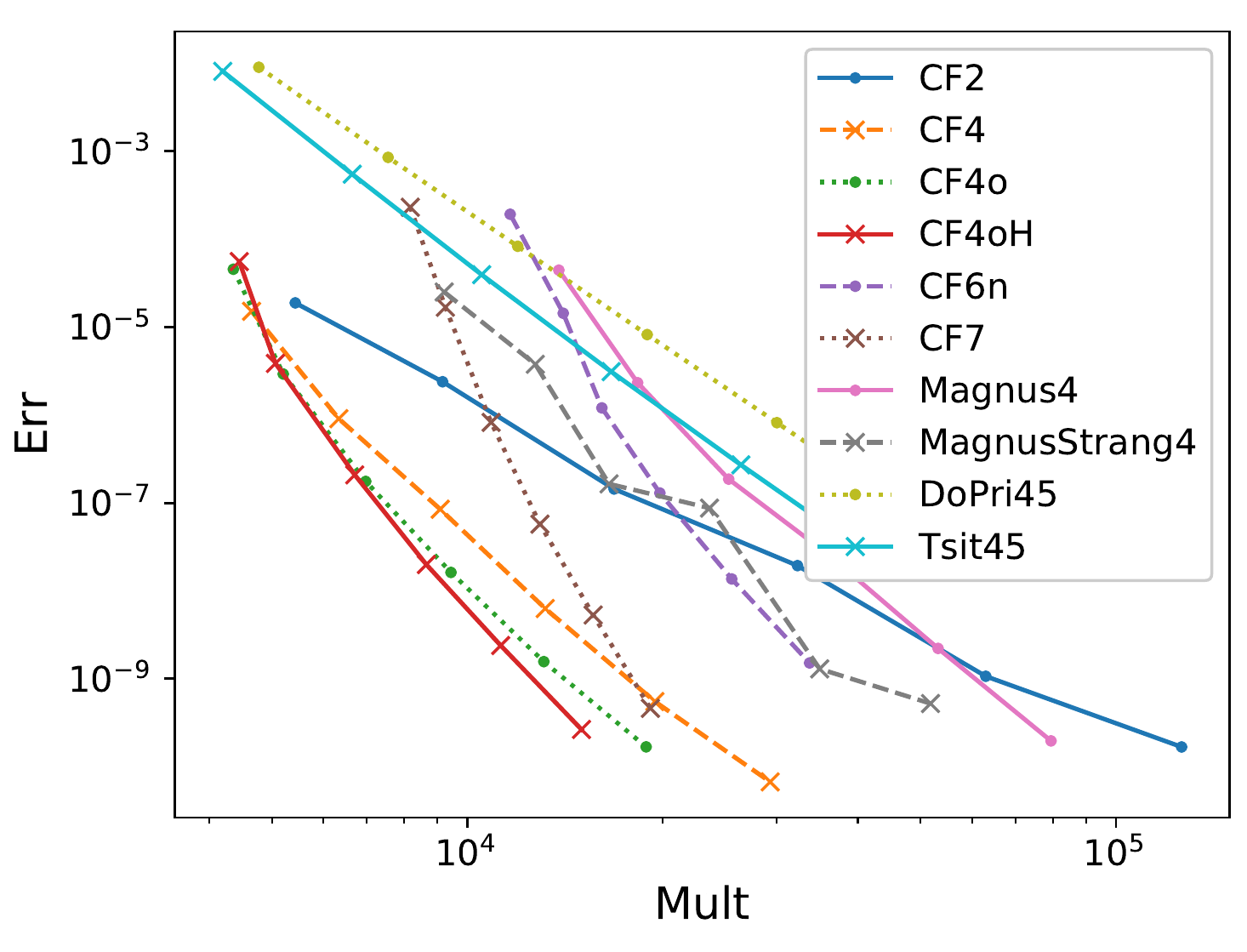} \quad \includegraphics[width=5.8cm]{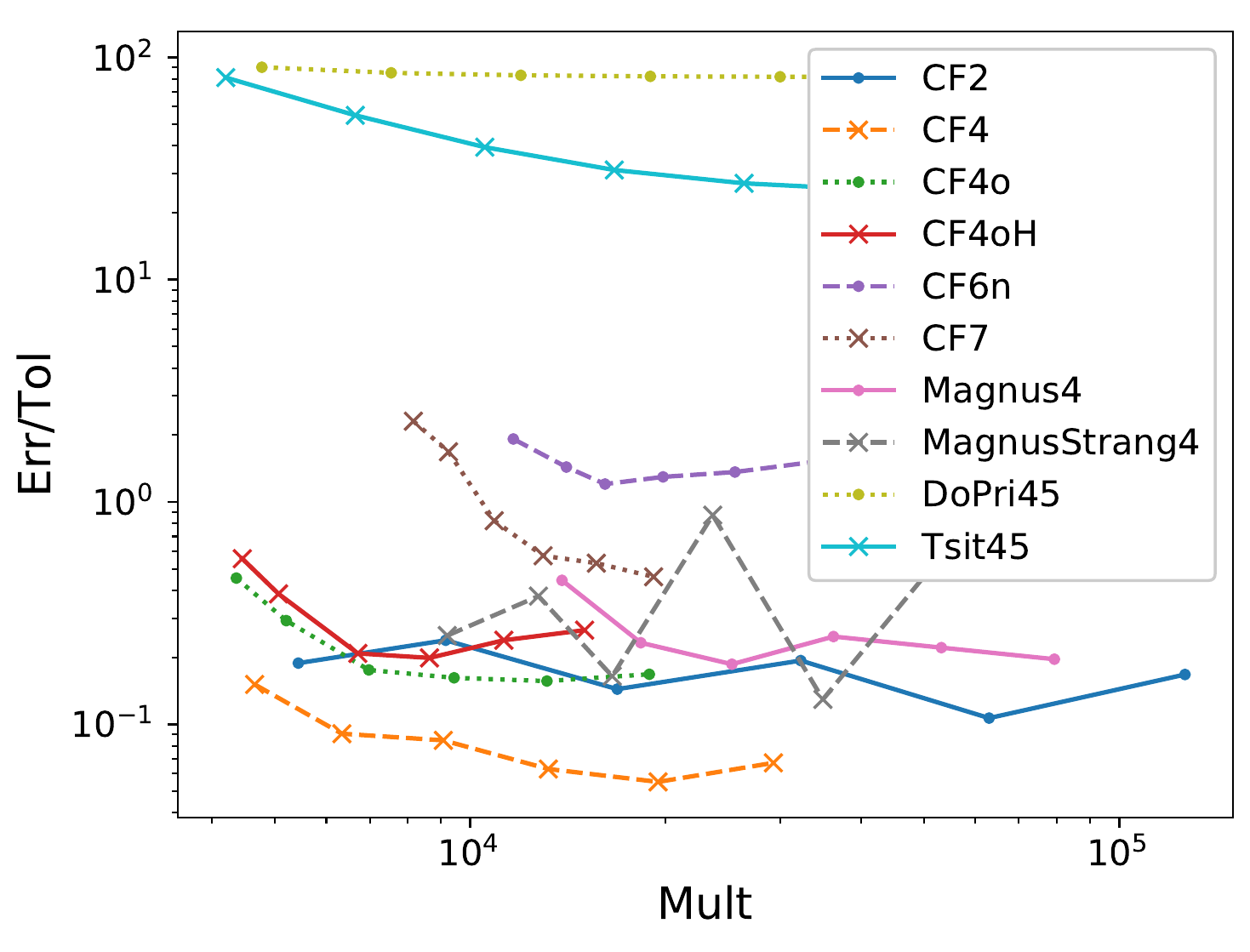}
\caption{$2\times4$ geometry, adaptive time-stepping. Error as a function of matrix-vector
multiplications (left) and quotient $\frac{\mathrm{error}}{\mathrm{tolerance}}$ (right).
\label{adapt8x1}}
\end{center}
\end{figure}

\begin{figure}[ht!]
\begin{center}
\includegraphics[width=5.8cm]{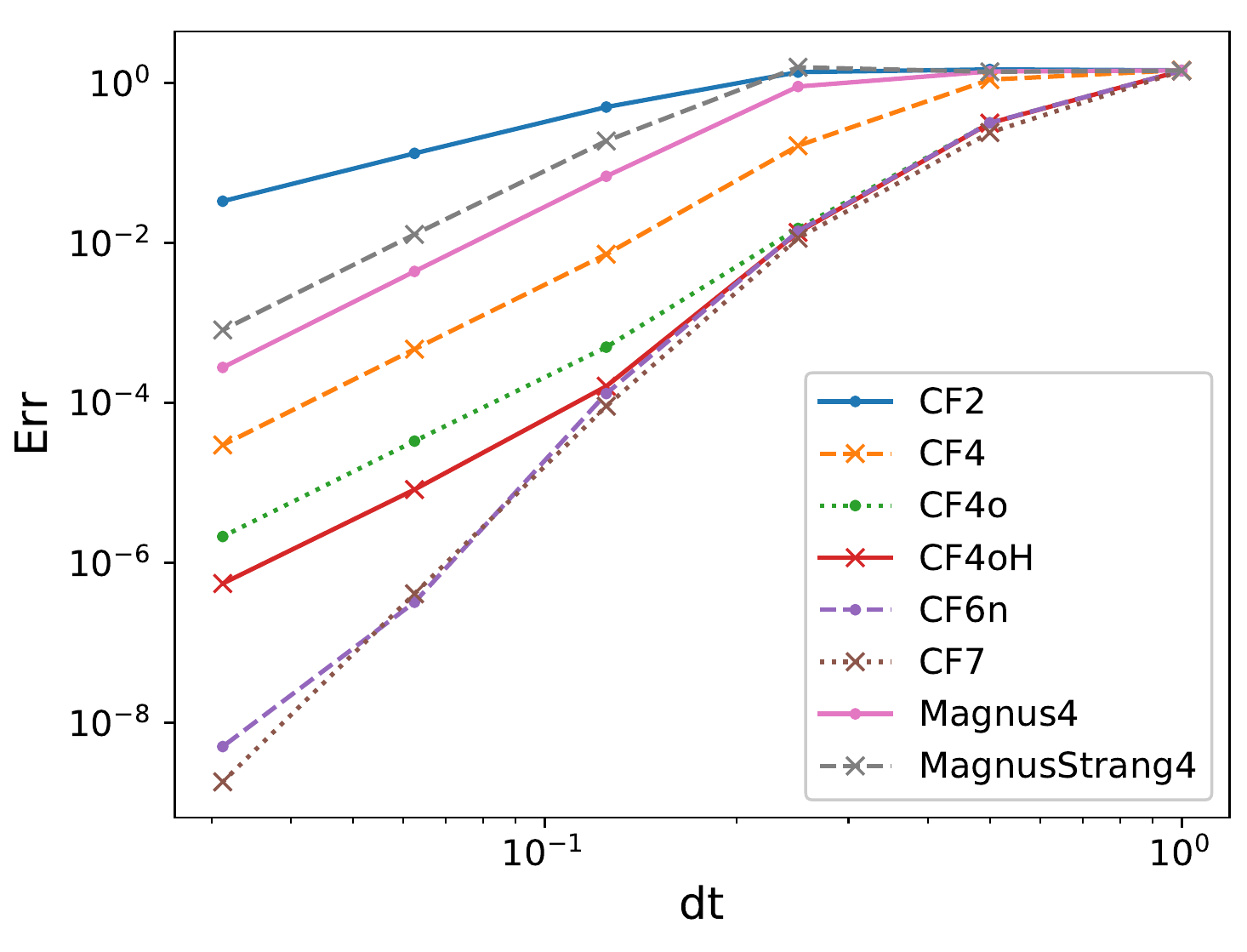} \quad \includegraphics[width=5.8cm]{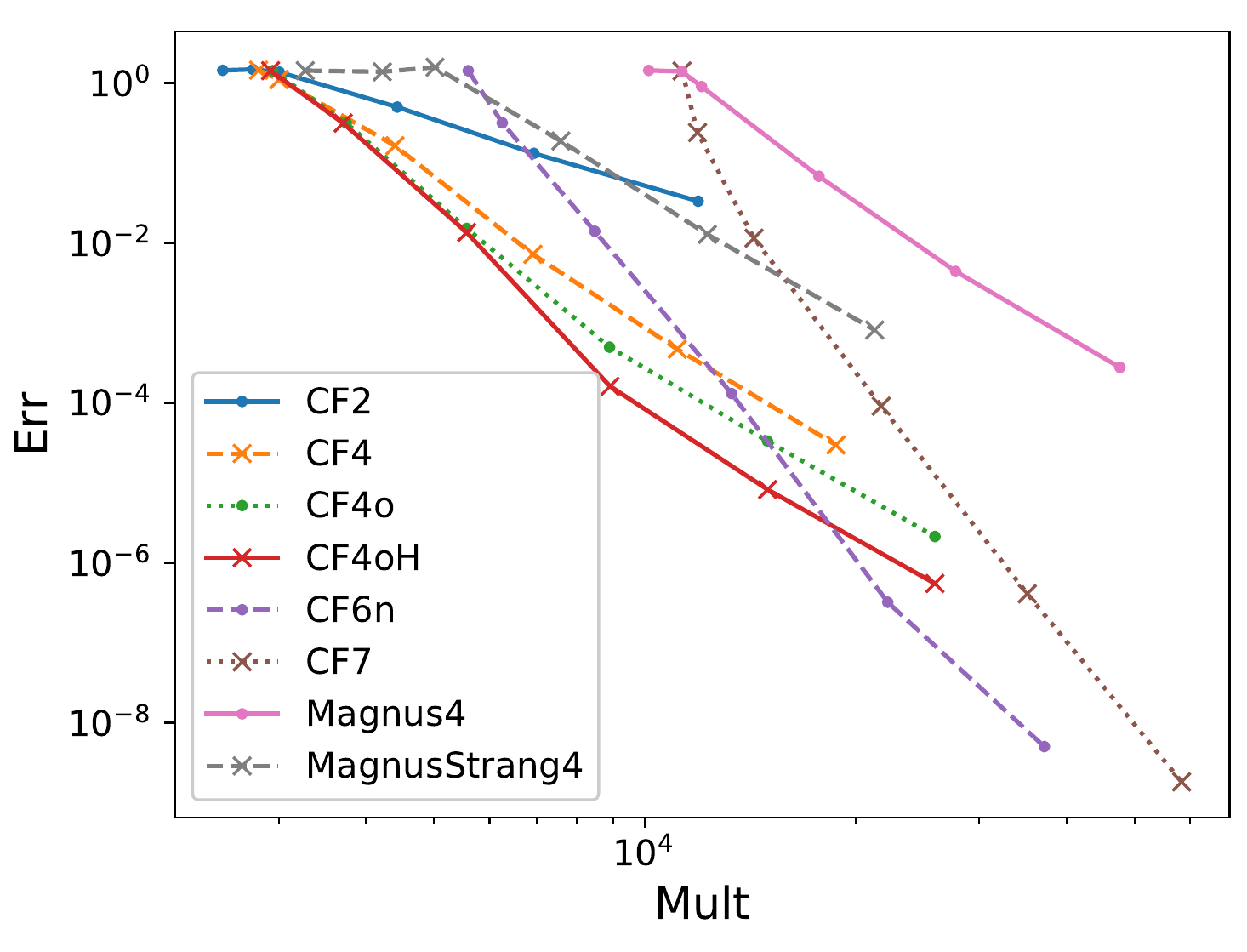}
\caption{ $4\times 3$ geometry, equidistant time-steps. Error as a function of the step-size (left) and as
a function of matrix-vector multiplications (right). \label{equi4x3}}
\end{center}
\end{figure}

\begin{figure}[ht!]
\begin{center}
\includegraphics[width=5.8cm]{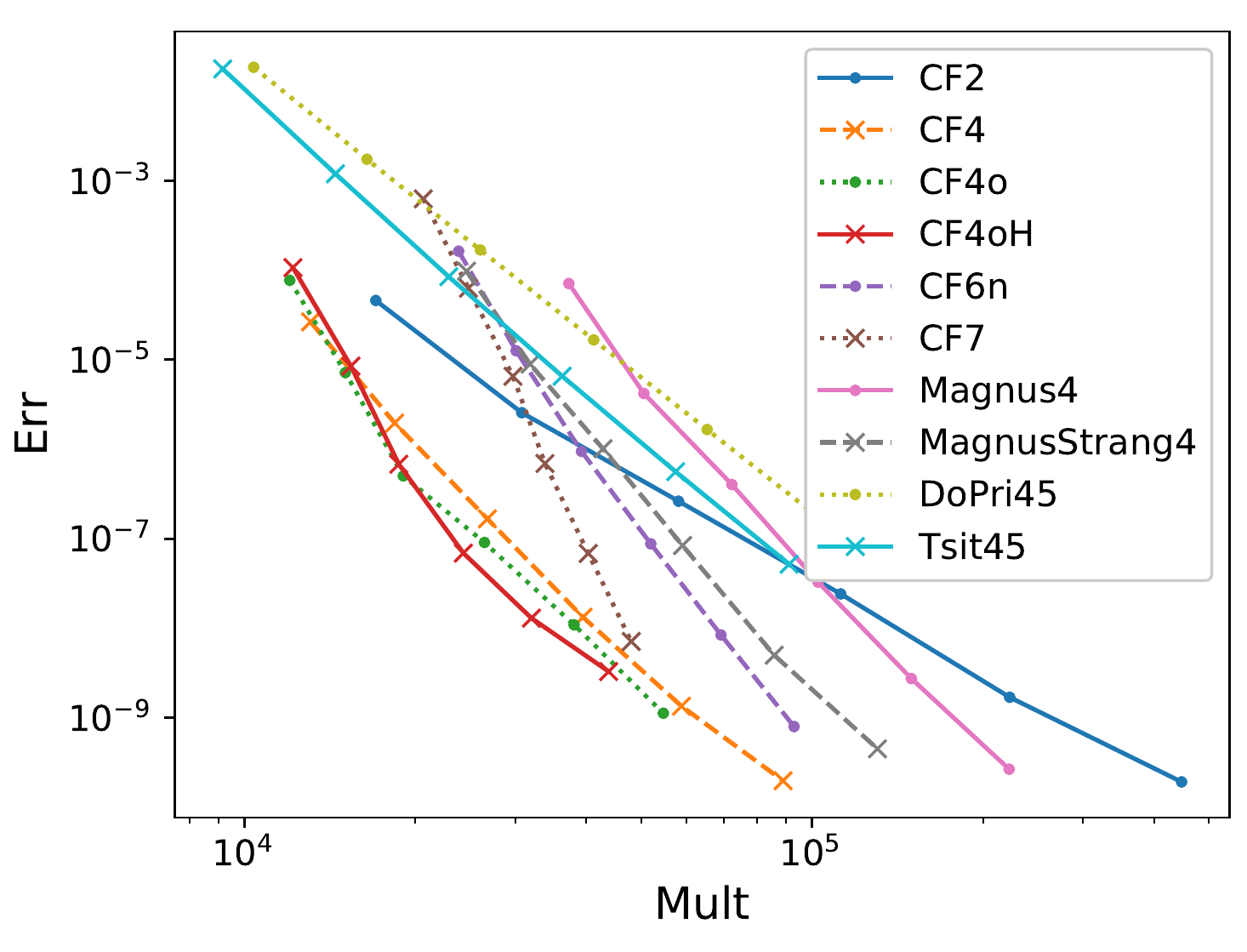} \quad \includegraphics[width=5.8cm]{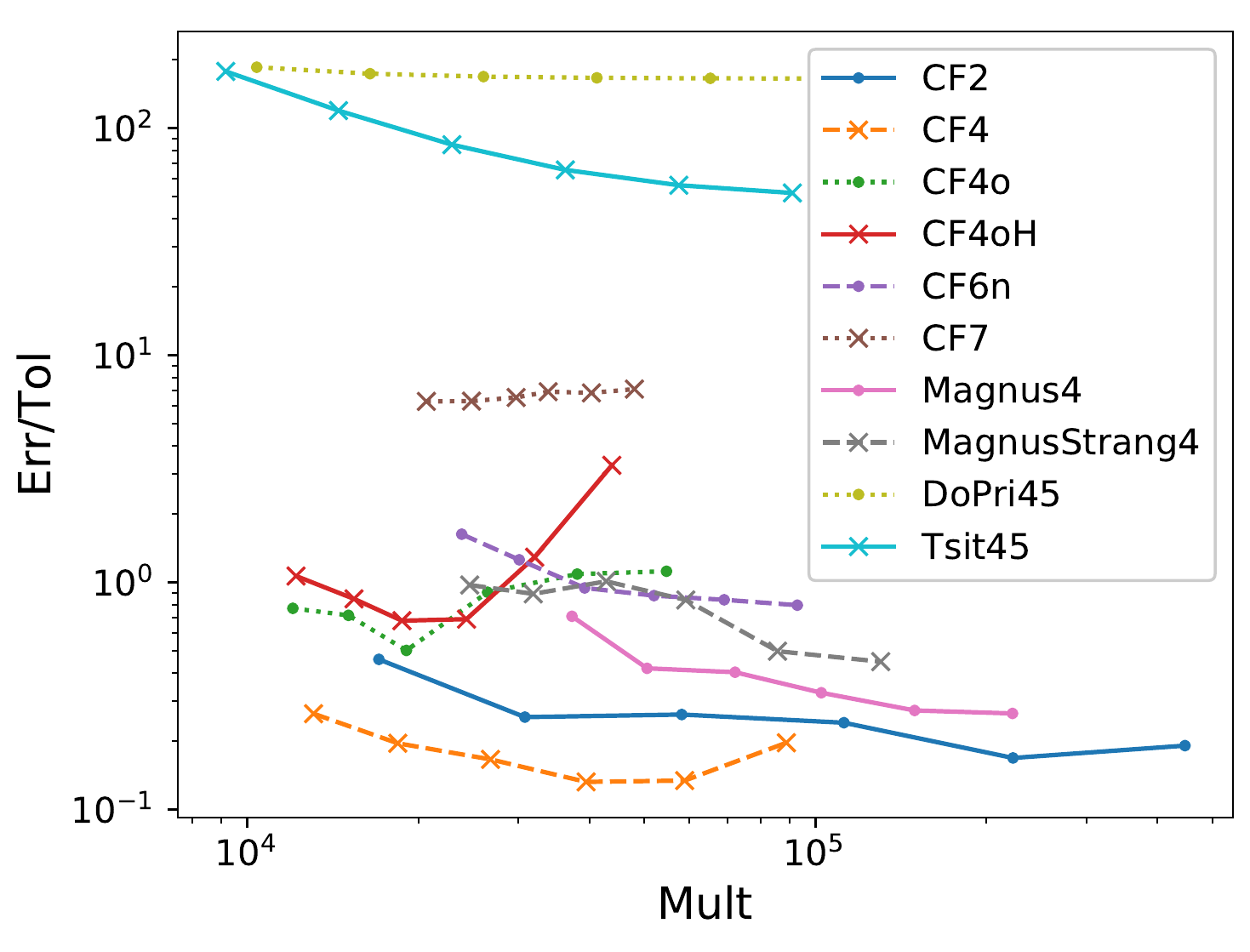}
\caption{$4\times 3$ geometry, adaptive time-stepping. Error as a function of matrix-vector
multiplications (left) and quotient $\frac{\mathrm{error}}{\mathrm{tolerance}}$ (right).
\label{adapt4x3}}
\end{center}
\end{figure}

First we study the $2\times 4$ ladder geometry in Figures~\ref{equi8x1} and~\ref{adapt8x1}, second a  $4\times 3$ lattice in
Figures~\ref{equi4x3} and~\ref{adapt4x3}.

Figures~\ref{equi8x1} and~\ref{equi4x3} show in the left plots the error of the
computation as a function of equidistant time-steps. The step-sizes are
consecutively halved such that $\tau= 2^{-k},\ k=0,\dots, 5$. The error here is
computed relative to a reference solution which was obtained with a tolerance $10^{-11}$. We observe that the theoretical
convergence orders are well reflected in the empirically determined orders
from the numerical experiments. Note the apparent order reduction for the
highest precisions, which is to be attributed to the limited accuracy
of the reference solution.

The right plots in Figures~\ref{equi8x1} and~\ref{equi4x3} show the error
as a function of matrix-vector multiplications in the computation
of the matrix exponentials. In fact, most of the computational effort
arising in the context of exponential integrators is associated with
this source. Observe that indeed, the plots show a systematic convergence
order similarly as the error/$dt$ plots. The results are, however,
also influenced by the fact that the Lanczos method requires more
iterations when the time-step is larger, see for example~\cite{hochbrucketal98,lubich07,niewri12,saad92s}.

Generally, the highest-order methods \texttt{CF7} and \texttt{CF6n}
have the highest accuracy, where the benefits of the order 7 approximation
are still superimposed by a larger error constant and \texttt{CF6n}
fares better in terms of the number of matrix-vector multiplications.
The advantages of the high-order methods are manifest especially for
high accuracies. A comparison of the fourth-order methods show that
\texttt{CF4o} and \texttt{CF4oH} are almost indistinguishable for
the $2\times4$ geometry, but the latter is slightly more accurate for the
$4\times3$ model in Fig.\ref{equi4x3}. Both the classical \texttt{Magnus4} integrator and
the new \texttt{MagnusStrang4} method suffer from large error
constants. The envisaged advantage manifests itself in the accuracy
as compared to the number of matrix-vector multiplications, however.
The new integrator thus shows its advantage over the classical
Magnus integrator, as was proposed in \cite{karo18a}.
Particularly it reduces the effort for the Lanczos process, but
the optimized commutator-free methods are clearly advantageous,
especially in their high-order variants. The second order exponential
midpoint rule is not competitive.

In this study, we put an emphasis on the adaptive implementation
of the time-stepping methods. Thus, we compare the achieved accuracy
with the number of matrix-vector multiplications in the respective
left-hand plots in Figures~\ref{adapt8x1} and~\ref{adapt4x3}.
We observe that apparently, the error estimators for the high-order
methods are expensive to compute, and the fourth order methods are
most efficient if adaptive time-stepping is included. For the reason of comparisons, we have also
tested adaptive Runge--Kutta methods which are very popular as
they are easy to implement and use, and state-of-the-art-implementations
are widely available. However, for the problem class under consideration,
both the Dormand/Prince method \cite{dopri80} (\texttt{DoPri45}
in the graphics) and an improved method by
Tsitouras \cite{tsitouras11} \texttt{Tsit45}
fail significantly to achieve the
prescribed tolerance, and are thus the least reliable,
and also inefficient integrators.

For this comparison of the usefulness of adaptive strategies,
reliability is also an important criterion. In the right-hand side plots
in Figures~\ref{adapt8x1} and~\ref{adapt4x3}, we thus show the
quotients of the achieved accuracy over the prescribed accuracy
(as a function of the number of matrix multiplications).
Values significantly larger than 1 indicate a very unreliable adaptive strategy
which fails to reach the prescribed tolerance, while values smaller than
1 signify an inefficient procedure which induces more computational
effort than is actually required to reach the tolerance. We observe that
for most exponential-based methods, this quotient is slightly smaller than 1,
where the highest-order methods are the most unreliable, but the
adaptive Runge--Kutta methods fail by far to satisfy the tolerance
criterion. The non-optimized commutator-free Magnus-type integrators \texttt{CF4}
and (to a lesser degree) \texttt{CF2} fall below the prescribed
tolerance requirement most noticeably.

From the experiments, we can rule out explicit Runge--Kutta
methods as appropriate integrators for the problem class we consider.
The picture is similar for both the tested models.

Figures~\ref{equi8x1}--\ref{adapt4x3} also allow an assessment
of the computational advantage of adaptive time-stepping over
equidistant grids. If we compare equal levels on the $y$-axis
in Figures~\ref{equi8x1} (right plot) and \ref{adapt8x1} (left plot),
for the $2\times 4$ geometry, and likewise Figures~\ref{equi4x3} (right plot)
and~\ref{adapt4x3}  (left plot) for the $4\times3$ geometry for the same integrators, we observe that
the same error level can be obtained at a smaller computational
effort in the adaptive computations for the majority of
time integrators. Additional tests, performed for different choices
of the parameters $\sigma_p$ and $\omega$ are given in \ref{further}.
These confirm the picture inferred here.

Finally, we illustrate the appropriateness of our step-size selection strategy
in the sense that the step-sizes are chosen according to the local smoothness
of the solution. We use the $2\times 4$ ladder geometry in this experiment.

\begin{figure}[ht!]
\begin{center}
\includegraphics[width=5.8cm]{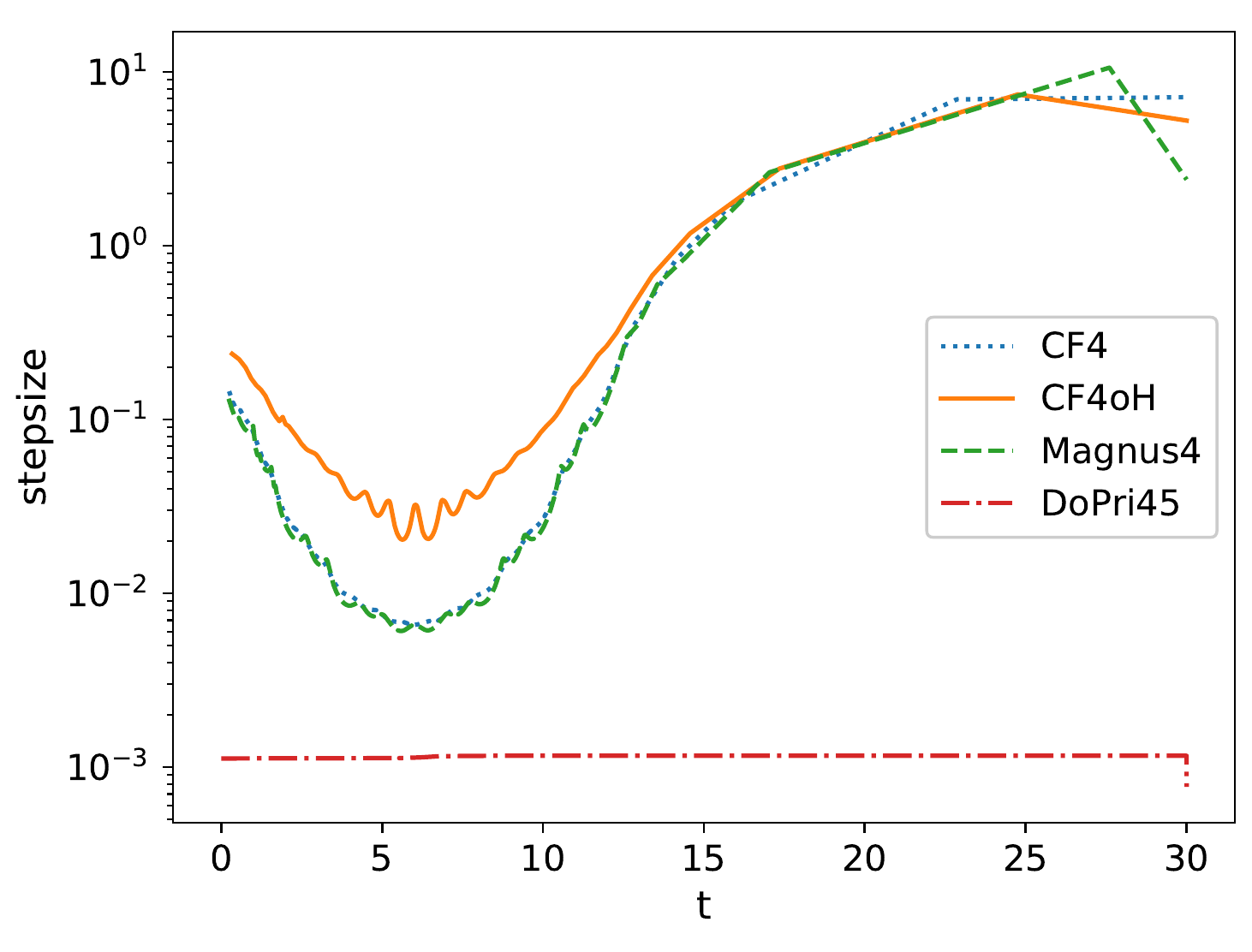}
\caption{$2\times4$ geometry. Adaptively chosen step-sizes for several integrators. \label{stepsizes}}
\end{center}
\end{figure}

Figure~\ref{stepsizes} shows the stepsizes chosen in the course of the time
integration. The local step-sizes \texttt{stepsizes} are shown as a
function of time $t$ for the representative choice of methods
\texttt{DoPri45},\ \texttt{CF4},\ \texttt{CF4oH}, and \texttt{Magnus4}.
We observe that the Runge--Kutta method chooses by far the smallest time-steps,
the optimized commutator-free Magnus-type method \texttt{CF4oH} allows
the largest time-steps and the classical Magnus integrator \texttt{Magnus4}
and \texttt{CF4} are comparable, with time-steps in between the other two methods.

\begin{figure}[ht!]
\begin{center}
\includegraphics[width=5.8cm]{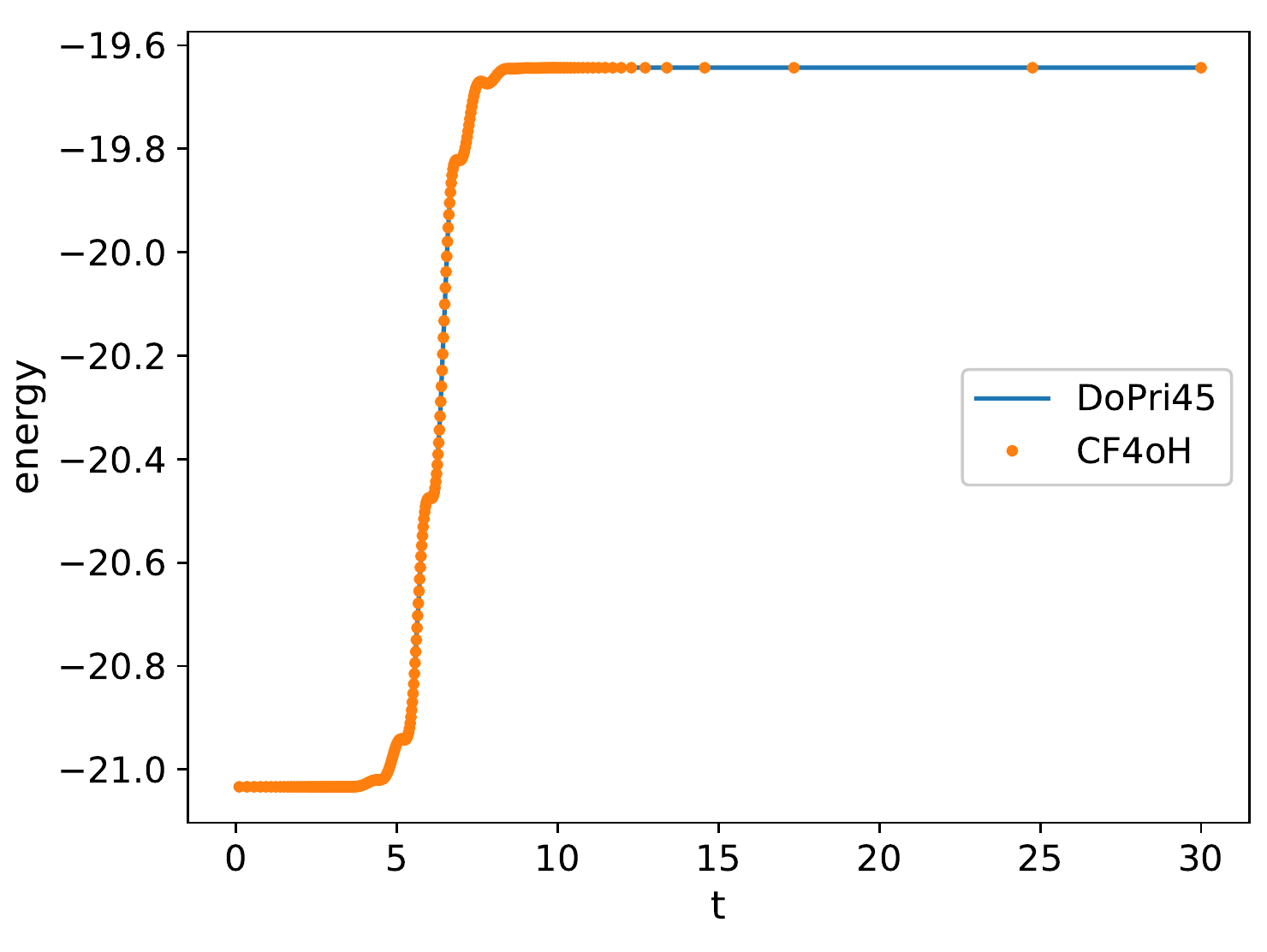} \quad \includegraphics[width=5.8cm]{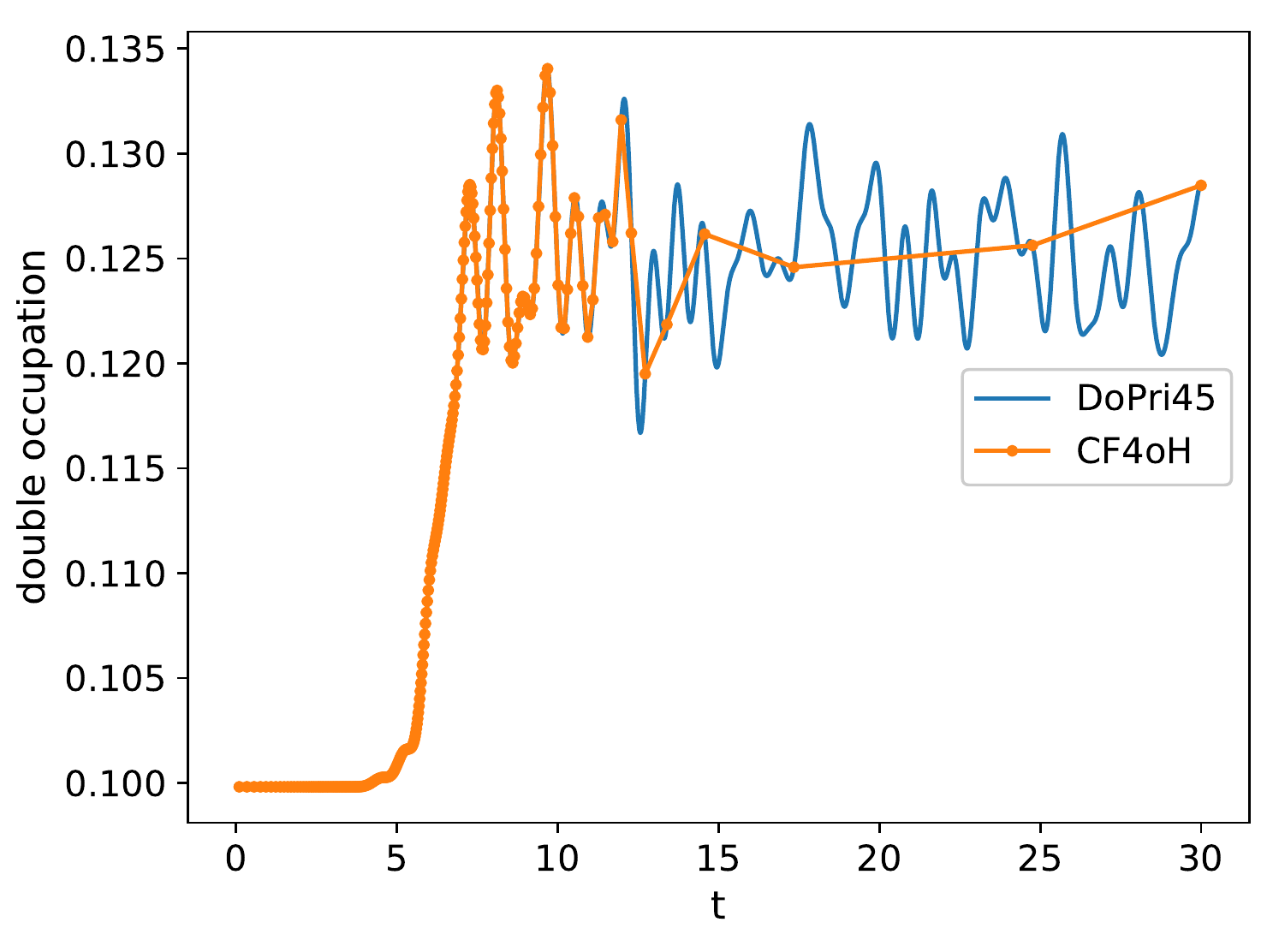}
\caption{$2\times4$ geometry. Approximation of the energy functional (left) and
the mean double occupation functional (right) by optimized fourth-order commutator free Magnus integrators (\texttt{CF4oH}) and RRunge-Kutta (\texttt{DoPri45}). \label{functionals}}
\end{center}
\end{figure}

Figure~\ref{functionals} shows the approximation quality of the two functionals
\emph{energy} and \emph{double occupation} by the commutator-free Magnus-type
method \texttt{CF4oH}. These quantities describe the energy transfer into the system and the number of electron-hole pairs (or double/single occupied sites) excited by the solar light, respectively.   For both the reference solution computed by \texttt{DoPri45}
and for \texttt{CF4oH} a tolerance of $10^{-11}$ was imposed.
The Runge--Kutta method chooses vanishingly small time-steps, thus the
numerical solution is plotted as a solid line, and the dots along these
curves represent the points chosen by the adaptive \texttt{CF4oH} method.
We observe that at the beginning of the time propagation the time-steps
are chosen as quite small, after attenuation of the external pulse, however,
time-steps increase impressively, while still the approximation quality of the
double occupation functional is remarkable: The approximations computed with
the large time-steps chosen for \texttt{CF4oH} lie on the curve provided
by the reference method, however not following the oscillations in between
the solution points. Overall, 244 points are needed for \texttt{CF4oH}, while \texttt{DoPri45}
requires 26015 steps.

\section{Conclusions}

In this study, we have investigated the successful application of adaptive
Magnus-type exponential integrator for Hubbard models with a time-dependent electric field, modeling e.g. the impact of a solar photon. Similar Hamiltonians and time-dependencies  will also be found in other situations where  strongly interacting electrons are driven out-of-equilibrium by an external field. For such problems, commutator-free Magnus-type methods were found to be preferable
over methods directly based on the Magnus expansion.

It was found that all methods show their expected convergence
orders on coherent equidistant grids. However, the methods based on the Magnus expansion have larger error
constants, where a newly proposed integrator improves on the
classical fourth order Magnus integrator, but commutator-free methods
are clearly to be favored. The use of high-order methods
grants high accuracy for a given computational effort, especially when
very precise solutions are sought.

We have also tested adaptive strategies based on asymptotically correct
estimators of the local time-stepping error. Fourth order methods
were found to be the most recommendable choice, where a new optimized method first
presented in this work performs best. The popular and easy to implement explicit Runge--Kutta
methods are not suitable for our problem class.

In addition to efficiency,
reliability is a major motivation to use adaptive time-stepping strategies.
Our tests revealed that commutator-free Magnus-type methods excel
in that the prescribed tolerance and the actually achieved error are
close (with high-order methods underestimating the necessary step-length
more pronouncedly), while the classical Magnus integrator is too
pessimistic (thus choosing unnecessarily small time-steps). Explicit
Runge--Kutta methods are prohibitively unreliable.

Adaptive commutator-free Magnus-type methods are thus concluded
to be the best choice for a reliable and efficient time integrators
of Hubbard models of solar cells, with the best results for
optimized fourth-order methods. This is also manifested by
observing that important functionals of the solution like
energy and mean double occupation are very well approximated
even for large time-steps.

\section*{Acknowledgements}

This work was supported by the Austrian Science Fund (FWF) [grant number P 30819-N32].
The computations have been conducted on the Vienna Scientific Cluster (VSC).
The work of K. Kropielnicka has been financed by The National Center of
Science (grant 2016/23/D/ST1/02061).

\begin{appendix}

\section{Implementation details for the fourth order Magnus-Strang splitting} \label{sec:adapt}

In this section, we give details of an efficient implementation of the new
integrator proposed in Section~\ref{subsec:MS} and the associated defect-based error estimator.

\subsection{Basic integrator}\label{sec:algpranav}

For an effective numerical scheme we have to approximate the integrals in  the definitions
of $\nA(\tau; t_0)$ and $\nB(\tau; t_0)$. We obtain
\begin{equation}\label{eq:numerical_scheme}
\nS(\tau; t_0) = {\rm e}^{\frac12 \tau^2\Phi_\nB(\tau;t_0)} {\rm e}^{\tau\Phi_\nA(\tau;t_0)}
{\rm e}^{\frac12 \tau^2\Phi_\nB(\tau;t_0)},
\end{equation}
where
\begin{align*}
\Phi_\nA(\tau;t_0)
     &= -\ii\,\Hdiag-\ii\,\hat{c}_1(\tau;t_0)\Hsymm+\hat{s}_1(\tau;t_0)\Hanti,\\
\Phi_\nB(\tau;t_0) &= -\hat{c}_2(\tau;t_0)[\Hsymm,\Hdiag]-\ii\,\hat{s}_2(\tau,t_0)[\Hanti,\Hdiag]-\ii\,\hat{r}(\tau;t_0)[\Hsymm,\Hanti]
\end{align*}
with
\begin{equation*}
\begin{split}
\hat{c}_1(\tau;t_0)&=\sum_{k=1}^{K_1}w^{(1)}_kc(t_0+x^{(1)}_k\tau)\approx\frac{1}{\tau}\tilde{c}_1(\tau,t_0)=\frac{1}{\tau}\int_{t_0}^{t_0+\tau}c(\zeta)\,{\rm d}\zeta,\\
\hat{s}_1(\tau;t_0)&=\sum_{k=1}^{K_1}w^{(1)}_kc(t_0+x^{(1)}_k\tau)\approx\frac{1}{\tau}\tilde{s}_1(\tau,t_0)=\frac{1}{\tau}\int_{t_0}^{t_0+\tau}s(\zeta)\,{\rm d}\zeta,\\
\hat{c}_2(\tau;t_0)&=\frac{1}{2}\sum_{k=1}^{K_2}w^{(2)}_k\big(c(t_0+y^{(2)}_k\tau)-c(t_0+x_k^{(2)}\tau)\big)\\
&\approx\frac{1}{\tau^2}\tilde{c}_2(\tau,t_0)=\frac{1}{2\tau^2}\int_{t_0}^{t_0+\tau}\int_{t_0}^{t_0+\zeta}\big(c(\zeta)-c(\xi)\big)\,{\rm d}\xi\,{\rm d}\zeta,\\
\hat{s}_2(\tau;t_0)&=\frac{1}{2}\sum_{k=1}^{K_2}w^{(2)}_k\big(s(t_0+y^{(2)}_k\tau)-s(t_0+x_k^{(2)}\tau)\big)\\
&\approx\frac{1}{\tau^2}\tilde{s}_2(\tau,t_0)=\frac{1}{2\tau^2}\int_{t_0}^{\tau}\int_{t_0}^{t_0+\zeta}\big(s(\zeta)-s(\xi)\big)\,{\rm d}\xi\,{\rm d}\zeta,\\
\hat{r}(\tau;t_0)&=\frac{1}{2}\sum_{k=1}^{K_2}w^{(2)}_k\big(c(t_0+y^{(2)}_k\tau)s(t_0+x^{(2)}_k\tau)-c(t_0+x_k^{(2)}\tau)s(t_0+y^{(2)}_k\tau)\big)\\
&\approx\frac{1}{\tau^2}\tilde{r}(\tau,t_0)=\frac{1}{2\tau^2}\int_{t_0}^{t_0+\tau}\int_{t_0}^{t_0+\zeta} \big(c(\zeta)s(\xi)-c(\xi)s(\zeta)\big)\,{\rm d}\xi\,{\rm d}\zeta, \\
\end{split}
\end{equation*}
where $x^{(1)}_k$, $w^{(1)}_k$ respectively $x^{(2)}_k$, $y^{(2)}_k$, $w^{(2)}_k$  are the nodes and weights of a suitable quadrature formula for integration over the interval $[0,1]$ and over the
triangle $0\leq y\leq 1$, $0\leq x\leq y$, respectively, see Tables~\ref{tab:quadcoef1}, \ref{tab:quadcoef2}.
\begin{table}
\begin{center}
\begin{tabular}{cc}
\hline
$x_k^{(1)}$ & $w_k^{(1)}$ \\
\hline
$\frac{1}{2}-\sqrt{\frac{1}{12}}$ & $\frac{1}{2}$ \\
$\frac{1}{2}+\sqrt{\frac{1}{12}}$ & $\frac{1}{2}$ \\
\hline
\end{tabular}\end{center}
\caption{Quadrature coefficients of order 4 for the interval $[0,1]$.\label{tab:quadcoef1}}
\end{table}
\begin{table}
\begin{center}
\begin{tabular}{ccc}
\hline
$x_k^{(2)}$ & $y_k^{(2)}$ & $w_k^{(2)}$\\
\hline
0.445948490915965 & 0.554051509084035 & 0.111690794839006\\
0.445948490915965 & 0.891896981831930 & 0.111690794839006 \\
0.108103018168070 & 0.554051509084035 & 0.111690794839006 \\
0.091576213509771 & 0.908423786490229 & 0.054975871827661 \\
0.091576213509771 & 0.183152427019541 & 0.054975871827661 \\
0.816847572980459 & 0.908423786490229 & 0.054975871827661 \\
\hline
\end{tabular}\end{center}
\caption{Quadrature coefficients of order 4 for the triangle $0\leq y\leq 1$, $0\leq x\leq y$.\label{tab:quadcoef2}}
\end{table}
Note that for later convenience (for the definition of a defect-based error estimator) we extracted
factors $\tau$ and $\tau^2$ in the exponents in (\ref{eq:numerical_scheme}).

The computation of the matrix--vector product
\begin{equation}\label{eq:yPhiv}
y=\Phi_\nB(\tau; t_0)\cdot v
\end{equation}
is accomplished using the algorithm in Table~\ref{lsmf}.
\begin{table}[ht]
\begin{center}
\begin{boxedminipage}{5cm}
\begin{tabbing}
$h_1 = \Hdiag\cdot v$\\
$h_2 = \Hsymm\cdot v$\\
$h_3 = \Hanti\cdot v$\\
$h_4 = -\ii\hat{c}_2h_1+\hat{r}h_3$\\
$y = \Hsymm\cdot h_4$\\
$h_4 = \hat{s}_2h_1-\hat{r}h_2$\\
$h_1 = \Hanti\cdot h_4$\\
$y = y + h_1$\\
$h_4 = -\ii\hat{c}_2h_2+\hat{s}_2h_3$\\
$h_1 = \Hdiag\cdot h_4 $\\
$y = -\ii(y-h_1)$
\end{tabbing}
\end{boxedminipage}
\caption{Algorithm for computing the matrix--vector product $y=\Phi_\nB(\tau; t_0)\cdot v$.\label{lsmf}}
\end{center}
\end{table}
In addition to cheap scaling, addition, and subtraction operations on vectors this algorithm incorporates
two (cheap) multiplications of vectors with $\Hdiag$ and two multiplications of vectors each with $\Hsymm$
and $\Hanti$. It needs four vectors $h_1,\dots, h_4$ for the storage of intermediate results.

\subsection{Implementation details for the symmetrized defect-based error estimator}\label{sec:implsymmdef}

The symmetrized defect here has the form
\begin{equation*}
\nD(\tau) = (\partial_\tau-\tfrac{1}{2}\partial_{t_0})\nS(\tau; t_0)+\tfrac{1}{2}\ii\big(
\Hfull(t_0+\tau)\nS(\tau; t_0)+\nS(\tau; t_0)\Hfull(t_0)
\big),
\end{equation*}
where $\partial_\tau$ denotes the derivative with respect to the first argument $\tau$, and
$\partial_{t_0}$ denotes the derivative with respect to the second argument $t_0$.
Here,
\begin{equation*}
\begin{split}
  (\partial_\tau-\tfrac{1}{2}\partial_{t_0})\nS(\tau; t_0) = &\Gamma_{\nB}(\tau;t_0){\rm e}^{\frac12\tau^2\Phi_\nB(\tau;t_0)}
  {\rm e}^{\tau\Phi_\nA(\tau,t_0)} {\rm e}^{\frac12 \tau^2\Phi_\nB(\tau;t_0)} \\
&+  {\rm e}^{\frac12 \tau^2\Phi_\nB(\tau;t_0)} \Gamma_\nA(\tau;t_0){\rm e}^{\tau\Phi_\nA(\tau;t_0)} {\rm e}^{\frac12 \tau^2\Phi_\nB(\tau;t_0)}\\
&+ {\rm e}^{\frac12 \tau^2\Phi_\nB(\tau;t_0)} {\rm e}^{\tau\Phi_\nA(\tau;t_0)} \Gamma_{\nB}(\tau;t_0){\rm e}^{\frac12\tau^2\Phi_\nB(\tau;t_0)}
\end{split}
\end{equation*}
with
\begin{align}
\Gamma_\nA(\tau)
&= \int_0^1{\rm e}^{\sigma\tau\Phi_\nA(\tau; t_0)} (\partial_\tau-\tfrac{1}{2}\partial_{t_0})\big(\tau\Phi_\nA(\tau;t_0)\big)
{\rm e}^{-\sigma\tau\Phi_\nA(\tau; t_0)}\,{\rm d}\sigma\nonumber\\
&=\Phi_\nA(\tau; t)+
\tau\int_0^1{\rm e}^{\sigma\tau\Phi_\nA(\tau; t_0)}\check{\Phi}_\nA(\tau; t_0){\rm e}^{-\sigma\tau\Phi_\nA(\tau; t_0)}\,{\rm d}\sigma,\label{eq:GammaA}
\end{align}
where
\begin{align*}
\check{\Phi}_\nA(\tau;t_0)&=(\partial_1-\tfrac{1}{2}\partial_2)\Phi_\nA(\tau; t_0)
= -\ii\,\check{c}_1(\tau;t_0)\Hsymm+\check{s}_1(\tau;t_0)\Hanti,\\
\check{c}_1(\tau;t_0)&=\sum_{k=1}^{K_1}w^{(1)}_k(x_k^{(1)}-\tfrac{1}{2})c'(t_0+x^{(1)}_k\tau),\\
\check{s}_1(\tau;t_0)&=\sum_{k=1}^{K_1}w^{(1)}_k(x_k^{(1)}-\tfrac{1}{2})s'(t_0+x^{(1)}_k\tau),
\end{align*}
and similarly,
\begin{align}
\Gamma_\nB(\tau; t_0)&=\int_0^1{\rm e}^{\sigma\frac{1}{2}\tau^2\Phi_\nB(\tau; t_0)}(\partial_\tau-\tfrac{1}{2}\partial_{t_0})
\big(\tfrac{1}{2}\tau^2\Phi_\nB(\tau;t_0)\big){\rm e}^{-\sigma\frac{1}{2}\tau^2\Phi_\nB(\tau; t_0)}\,{\rm d}\sigma
\nonumber \\
&= \tau\Phi_\nB(\tau; t_0)
+\frac{\tau^2}{2}\int_0^1{\rm e}^{\sigma\frac{1}{2}\tau^2\Phi_\nB(\tau; t_0)}
\check{\Phi}_\nB(\tau; t_0)
{\rm e}^{-\sigma\frac{1}{2}\tau^2\Phi_\nB(\tau; t_0)}\,{\rm d}\sigma\label{eq:GammaB}
\end{align}
with
\begin{equation*}
\check{\Phi}_\nB(\tau; t_0)
= -\check{c}_2(\tau;t_0)[\Hsymm,\Hdiag]-\ii\,\check{s}_2(\tau,t_0)[\Hanti,\Hdiag]-\ii\,\check{r}(\tau,t_0)[\Hsymm,\Hanti],
\end{equation*}
where
\begin{align*}
\check{c}_2(\tau;t_0)&=\frac{1}{2}\sum_{k=1}^{K_2}w^{(2)}_k
\big((y_k^{(2)}-\tfrac{1}{2})c'(t_0+y^{(2)}_k\tau)-(x_k^{(2)}-\tfrac{1}{2})c'(t_0+x_k^{(2)}\tau)
\big),\\
\check{s}_2(\tau,t_0)&=\frac{1}{2}\sum_{k=1}^{K_2}w^{(2)}_k
\big((y_k^{(2)}-\tfrac{1}{2})s'(t_0+y^{(2)}_k\tau)-(x_k^{(2)}-\tfrac{1}{2})s'(t_0+x_k^{(2)}\tau)
\big),\\
\check{r}(\tau,t_0)&=\frac{1}{2}\sum_{k=1}^{K_2}w^{(2)}_k
\big((y_k^{(2)}-\tfrac{1}{2})c'(t_0+y^{(2)}_k\tau)s(t_0+y^{(2)}_k\tau)
+(x_k^{(2)}-\tfrac{1}{2})c(t_0+y^{(2)}_k\tau)s'(t_0+y^{(2)}_k\tau)\\
&\qquad\qquad\ \ -(x_k^{(2)}-\tfrac{1}{2})c'(t_0+x^{(2)}_k\tau)s(t_0+y^{(2)}_k\tau)
-(y_k^{(2)}-\tfrac{1}{2})c(t_0+x^{(2)}_k\tau)s'(t_0+y^{(2)}_k\tau)\big).
\end{align*}

A computable approximation for (\ref{eq:GammaA}) is
\begin{equation*}
\widetilde{\Gamma}_\nA(\tau;t_0)=\Phi_\nA(\tau; t_0)+\sum_{m=0}^{3}\frac{\tau^{m+1}}{(m+1)!}
\mathrm{ad}_{\Phi_{\nA}(\tau;t_0)}^m(\check{\Phi}_\nA(\tau;t_0)) 
\end{equation*}
with leading error term
\begin{equation*}
\Gamma_\nA(\tau;t_0)-\widetilde{\Gamma}_\nA(\tau;t_0)=\frac{\tau^5}{5!}[\Phi_\nA,[\Phi_\nA,[\Phi_\nA,[\Phi_\nA,\check{\Phi}_\nA]]]+
\mathcal{O}(\tau^6).
\end{equation*}
From the order conditions for the quadrature coefficients
\begin{equation*}
\sum_{k=1}^{K_1}w_k^{(1)}=\int_0^1 1\,\mathrm{d}x=1,\quad\sum_{k=1}^{K_1}w_k^{(1)}x_k^{(1)}=\int_0^1 x\,\mathrm{d}x=\frac{1}{2}
\end{equation*}
it follows
\begin{equation*}
\check{c}_1(\tau,t_0)=\underbrace{\left(\sum_{k=1}^{K_1}w^{(1)}_k(x_k^{(1)}-\tfrac{1}{2})\right)}_{=0}c'(t_0)
+\mathcal{O}(\tau)=\mathcal{O}(\tau),
\end{equation*}
similarly $\check{s}_1(\tau,t_0)=\mathcal{O}(\tau)$, and thus
$\check{\Phi}_A(\tau;t_0)=\mathcal{O}(\tau)$.
We conclude
\begin{equation*}
\Gamma_\nA(\tau;t_0)-\widetilde{\Gamma}_\nA(\tau;t_0)=\mathcal{O}(\tau^6),
\end{equation*}
as required.

Similarly, a computable approximation for (\ref{eq:GammaB}) is
 \begin{equation*}
\widetilde{\Gamma}_\nB(\tau;t_0)=\tau\Phi_\nB(\tau;t_0)+\frac{\tau^2}{2}\check{\Phi}_\nB(\tau;t_0),
\end{equation*}
with leading error term
\begin{equation*}
\Gamma_\nB-\widetilde{\Gamma}_\nB=\frac{\tau^4}{8}[\Phi_\nB(\tau; t_0),\check{\Phi}_\nB(\tau;t_0)]+\mathcal{O}(\tau^6),
\end{equation*}
where we used $\Phi_\nB(\tau; t_0)=\mathcal{O}(\tau)$,   $\check{\Phi}_\nB(\tau, t_0)=\mathcal{O}(1)$. From the order conditions for the quadrature coefficients
\begin{align*}
&\sum_{k=1}^{K_2}w_k^{(2)}=\int_0^1\int_0^y 1\,\mathrm{d}x\,\mathrm{d}y=\frac{1}{2},\quad
\sum_{k=1}^{K_2}w_k^{(2)}x_k^{(2)}=\int_0^1\int_0^y x\,\mathrm{d}x\,\mathrm{d}y=\frac{1}{3},\\
&\sum_{k=1}^{K_2}w_k^{(2)}y_k^{(2)}=\int_0^1\int_0^y y\,\mathrm{d}x\,\mathrm{d}y=\frac{1}{3}
\end{align*}
it follows
\begin{equation*}
\check{c}_2(\tau;t_0)=\frac{1}{2}\underbrace{\left(\sum_{k=1}^{K_2}w^{(2)}_k
(y_k^{(2)}-x_k^{(2)})\right)}_{=0}c'(t_0)+\mathcal{O}(\tau)=\mathcal{O}(\tau),
\end{equation*}
similarly $\check{s}_2(\tau,t_0)=\mathcal{O}(\tau)$, $\check{r}(\tau,t_0)=\mathcal{O}(\tau)$, and thus
$\check{\Phi}_\nB(\tau;t_0)=\mathcal{O}(\tau)$.
We conclude
\begin{equation*}
\Gamma_\nB-\widetilde{\Gamma}_\nB=\mathcal{O}(\tau^6),
\end{equation*}
as required.
\begin{table}
\begin{center}
\begin{boxedminipage}{5cm}
\begin{tabbing}
        $h_2 = X\cdot v$\\
        $y = h_2$\\
        $h_1=Y\cdot v$\\
        $y = y + \tau h_1 $\\
        $h_3=X\cdot h_1$ /\!/ $=XYv$\\
        $y = y + \frac{1}{2}\tau^2 h_3 $\\
        $h_1=X\cdot h_3$ /\!/ $=XXYv$\\
        $y = y + \frac{1}{6}\tau^3 h_1 $\\
        $h_3=X\cdot h_1$ /\!/ $=XXXYv$\\
        $y = y + \frac{1}{24}\tau^4 h_3 $\\
        $h_1=Y\cdot h_2$ /\!/ $=YXv$\\
        $y = y - \frac{1}{2}\tau^2 h_1 $\\
        $h_3=X\cdot h_1$ /\!/ $=XYXv$\\
        $y = y - \frac{1}{3}\tau^3 h_3 $\\
        $h_1=X\cdot h_3$ /\!/ $=XXYXv$\\
        $y = y - \frac{1}{8}\tau^4 h_1 $\\
        $h_4=X\cdot h_2$ /\!/ $=XXv$\\
        $h_1=Y\cdot h_4$ /\!/ $=YXXv$\\
        $y = y + \frac{1}{6}\tau^3 h_1$\\
        $h_3=X\cdot h_1$ /\!/ $=XYXXv$\\
        $y = y + \frac{1}{8}\tau^4 h_3$\\
        $h_2=X\cdot h_4$ /\!/ $=XXXv$\\
        $h_1=Y\cdot h_2$ /\!/ $=YXXXv$\\
        $y = y - \frac{1}{24}\tau^4 h_1$
\end{tabbing}
\end{boxedminipage}\qquad
\begin{boxedminipage}{5cm}
\begin{tabbing}
$u=u_0$\\
$d=\frac{1}{2}\mathrm{i}\Hfull(t_0)\cdot u$\\
$u=\mathrm{e}^{\frac{1}{2}\tau^2\Phi_B(\tau;t_0)}u$\\
$d=\mathrm{e}^{\frac{1}{2}\tau^2\Phi_B(\tau;t_0)}d$\\
$d=d+\widetilde{\Gamma}_B(\tau;t_0)\cdot u$\\
$u=\mathrm{e}^{\tau\Phi_A(\tau;t_0)}u$\\
$d=\mathrm{e}^{\tau\Phi_A(\tau;t_0)}d$\\
$d=d+\widetilde{\Gamma}_A(\tau;t_0)\cdot u$\\
$u=\mathrm{e}^{\frac{1}{2}\tau^2\Phi_B(\tau;t_0)}u$\\
$d=\mathrm{e}^{\frac{1}{2}\tau^2\Phi_B(\tau;t_0)}d$\\
$d=d+\widetilde{\Gamma}_B(\tau;t_0)\cdot u$\\
$d=d+\frac{1}{2}\mathrm{i}\Hfull(t_0+\tau)\cdot u$
\end{tabbing}
\end{boxedminipage}
\end{center}
\caption{Left: algorithm for computing $y=\big(X+\sum_{m=0}^{3}\frac{\tau^{m+1}}{(m+1)!}
\mathrm{ad}_{X}^m(Y)\big)\cdot v$ which incorporates 13 matrix-vector multiplications with $X$ or $Y$ and
4 vectors $h_1,\dots, h_4$ for the storage of intermediate results.
 Right: algorithm for the simultaneous computation of the numerical solution $u=\mathcal{S}(\tau; t_0)u_0$ and the
 symmetrized defect $d=\mathcal{D}(\tau; t_0)u_0$. \label{tab:algorithms}}
\end{table}

Table~\ref{tab:algorithms} (right) shows the algorithmic realization of the symmetrized defect. The
applications of $\widetilde{\Gamma}_\nA(\tau;t_0)$ can be realized by the algorithm given in
Table~\ref{tab:algorithms} (left), the application of $\widetilde{\Gamma}_\nB(\tau;t_0)$ can be realized
by an obvious adaptation of the algorithm for (\ref{eq:yPhiv}).

\section{Further comparisons}\label{further}

To corroborate our conclusions about the different integration methods,
we give results showing the achieved accuracy as in Section~\ref{sec:numexp} but now for different choices of the parameters
$\sigma_p$ and $\omega$, i.e., for the length and frequency of the electric field pulse.

\begin{figure}[ht!]
\begin{center}
\includegraphics[width=5.8cm]{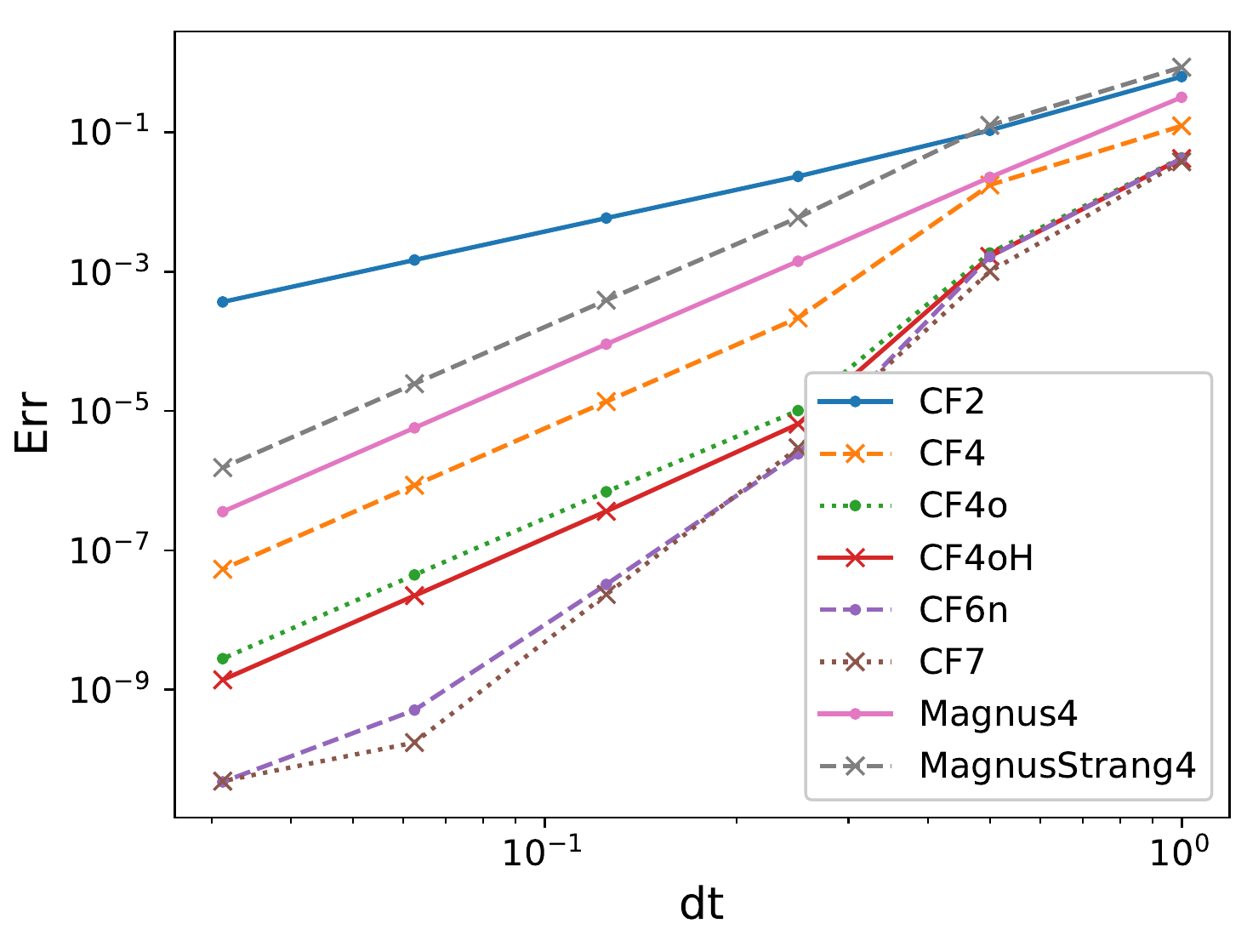} \quad \includegraphics[width=5.8cm]{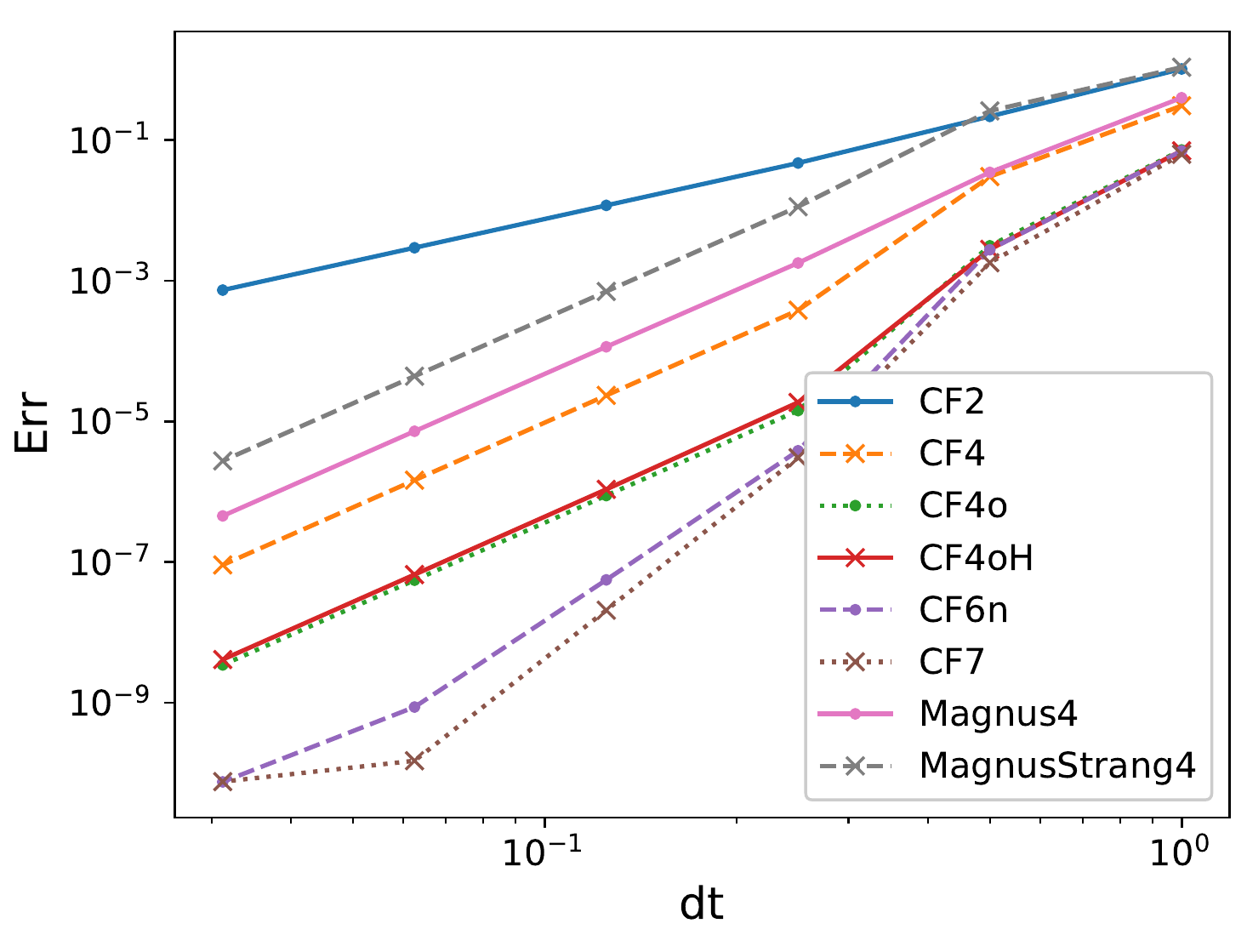} \\
\includegraphics[width=5.8cm]{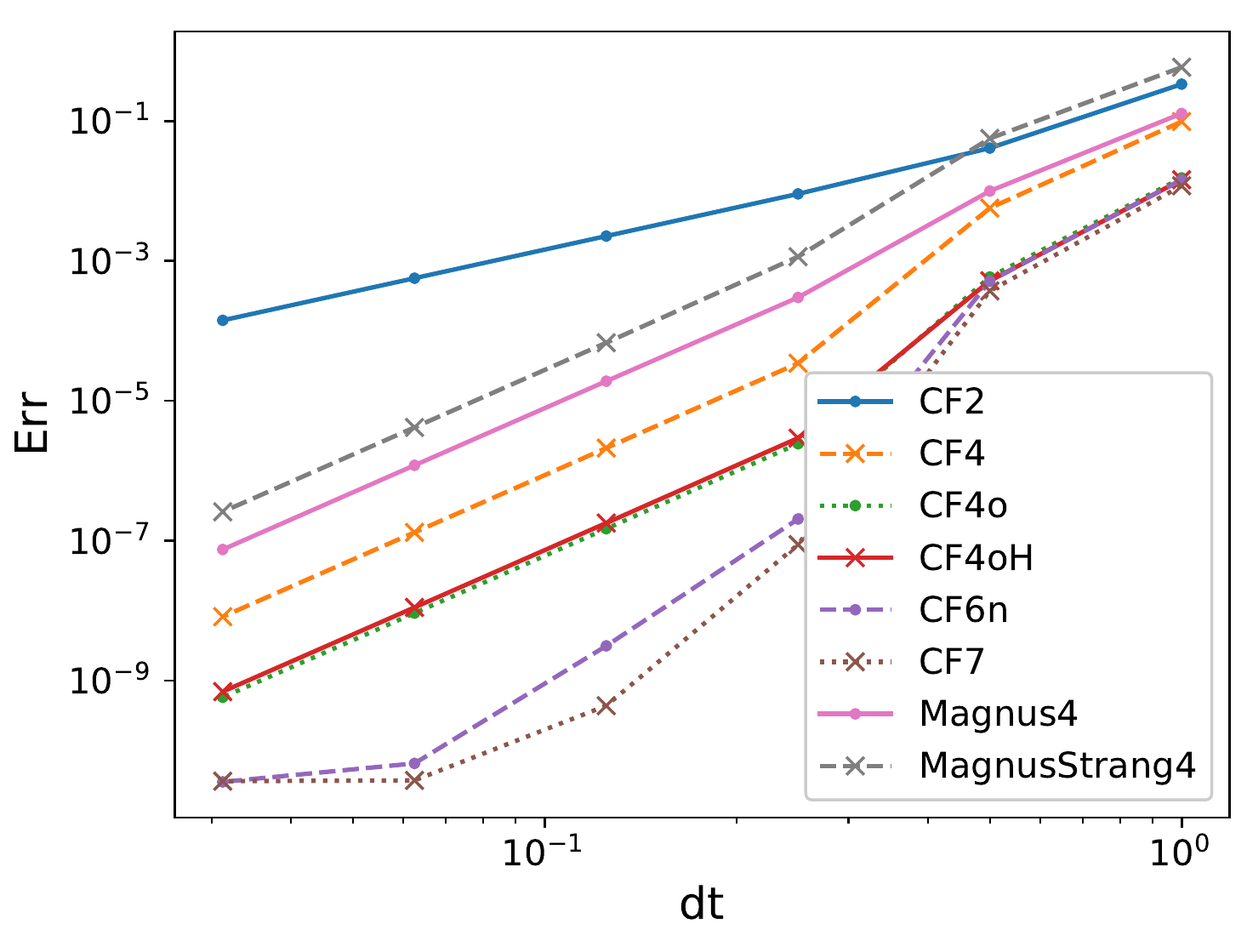} \quad \includegraphics[width=5.8cm]{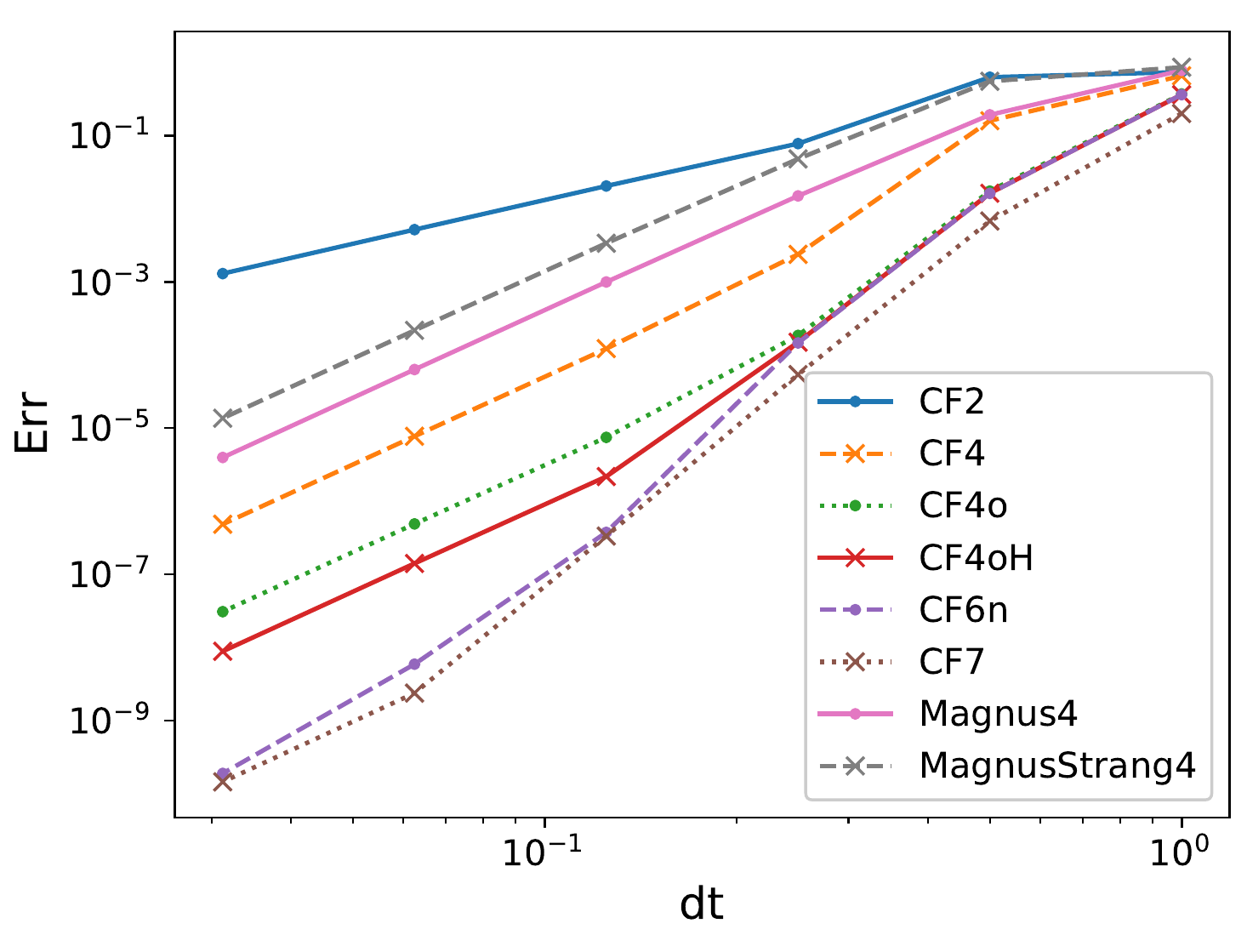} \\
\caption{$2\times4$ geometry, equidistant time-stepping. Error as a function of the step-size
for $\omega=3.5$ (top row), $\sigma_p=1$ (top left) and $\sigma_p=4$ (top right),
and for $\sigma_p=2$ (bottom row), $\omega=1.75$ (bottom left) and $\omega=7$ (bottom right).
\label{8x1variationequidist}}
\end{center}
\end{figure}

\begin{figure}[ht!]
\begin{center}
\includegraphics[width=5.8cm]{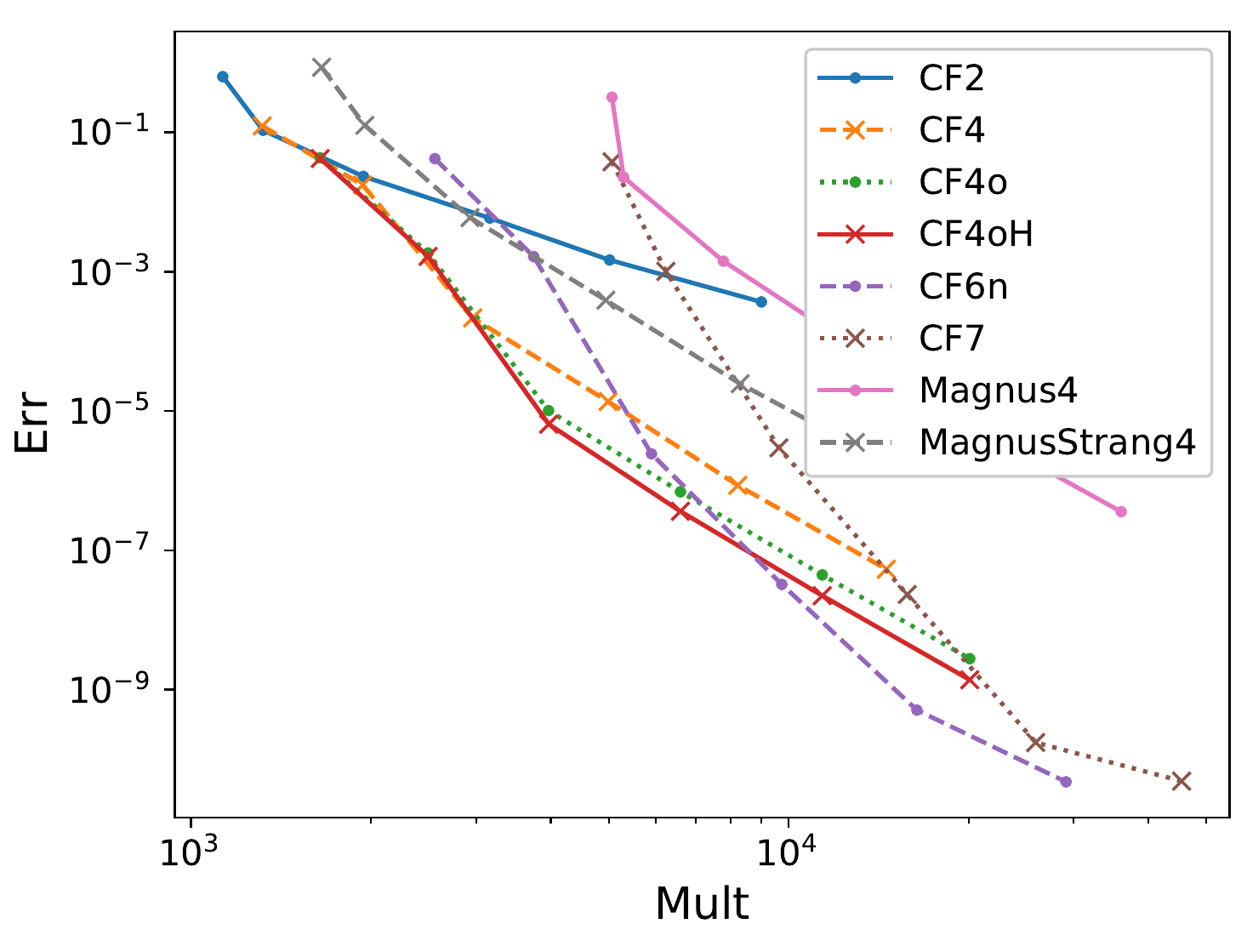} \quad \includegraphics[width=5.8cm]{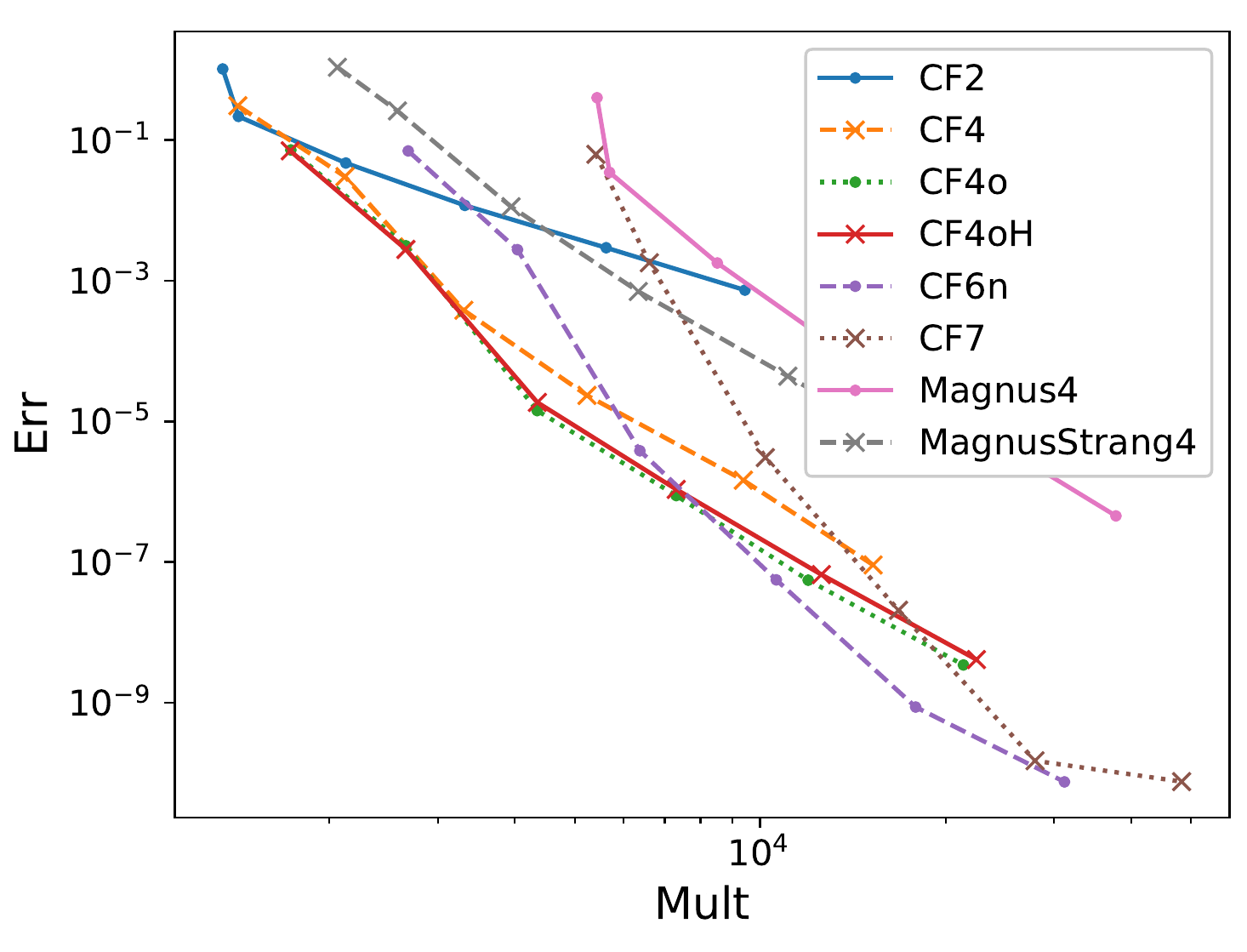} \\
\includegraphics[width=5.8cm]{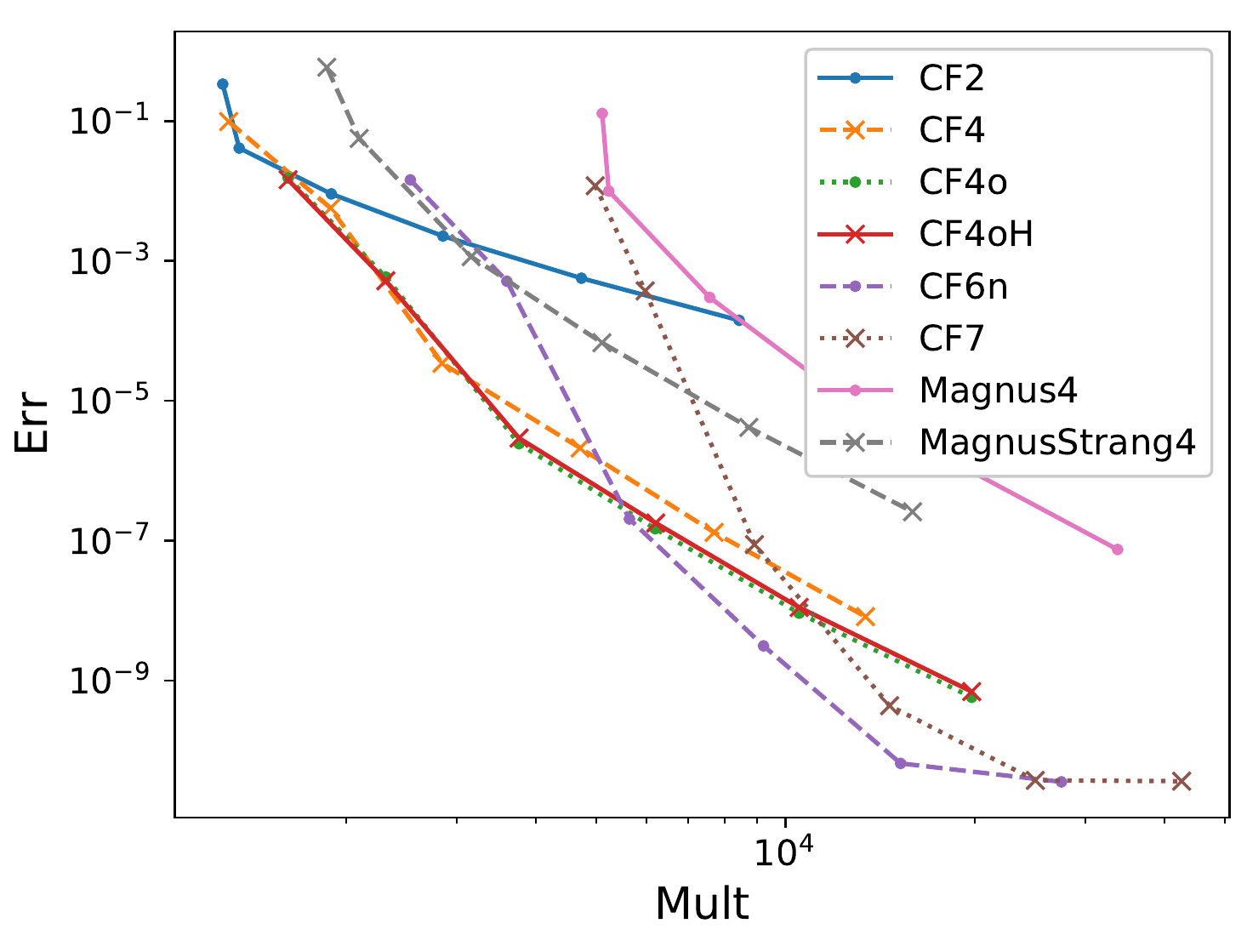} \quad \includegraphics[width=5.8cm]{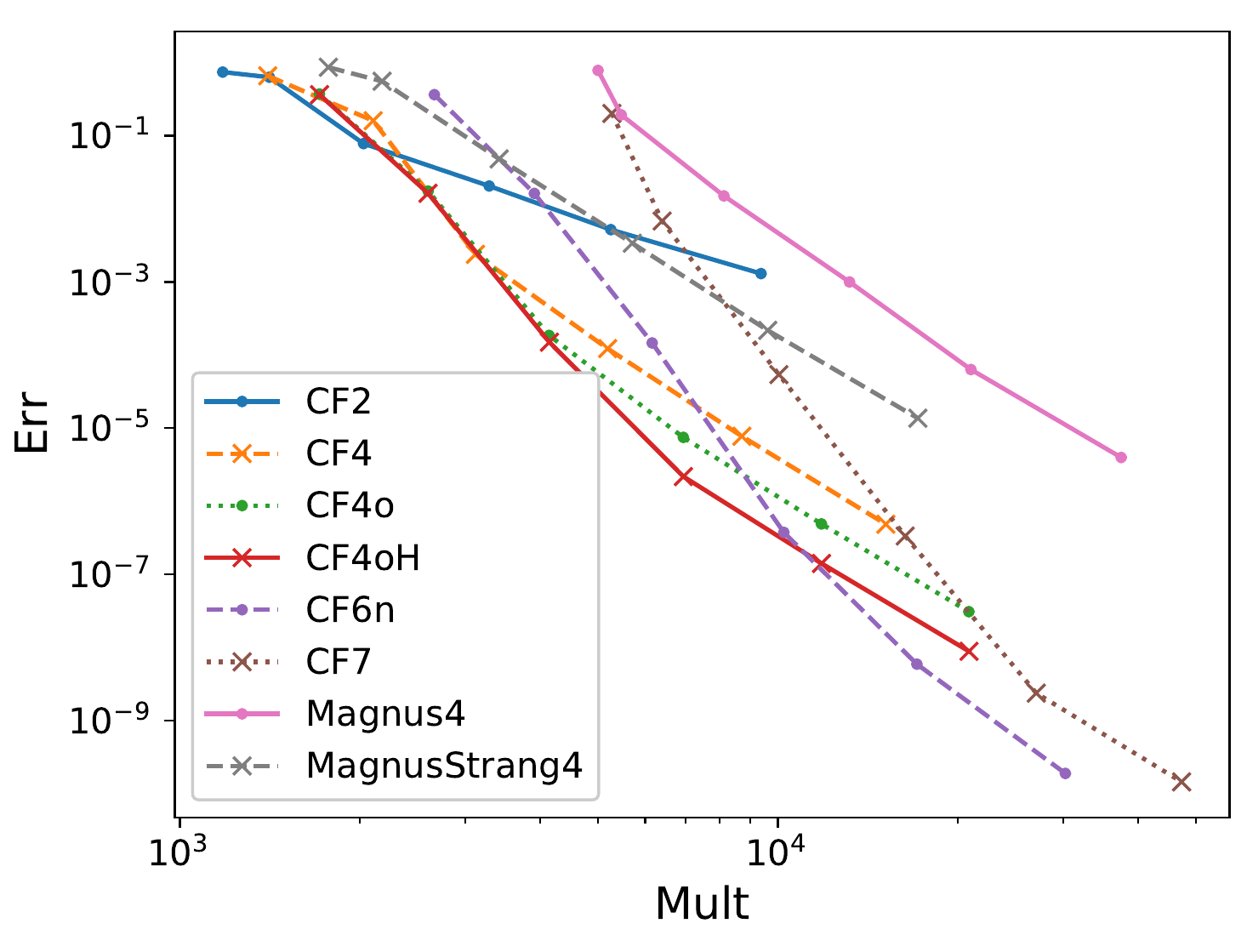} \\
\caption{$2\times4$ geometry, equidistant time-stepping. Error as a function of matrix--vector multiplications
for $\omega=3.5$ (top row), $\sigma_p=1$ (top left) and $\sigma_p=4$ (top right),
and for $\sigma_p=2$ (bottom row), $\omega=1.75$ (bottom left) and $\omega=7$ (bottom right).
\label{8x1variationequidistmult}}
\end{center}
\end{figure}

\begin{figure}[ht!]
\begin{center}
\includegraphics[width=5.8cm]{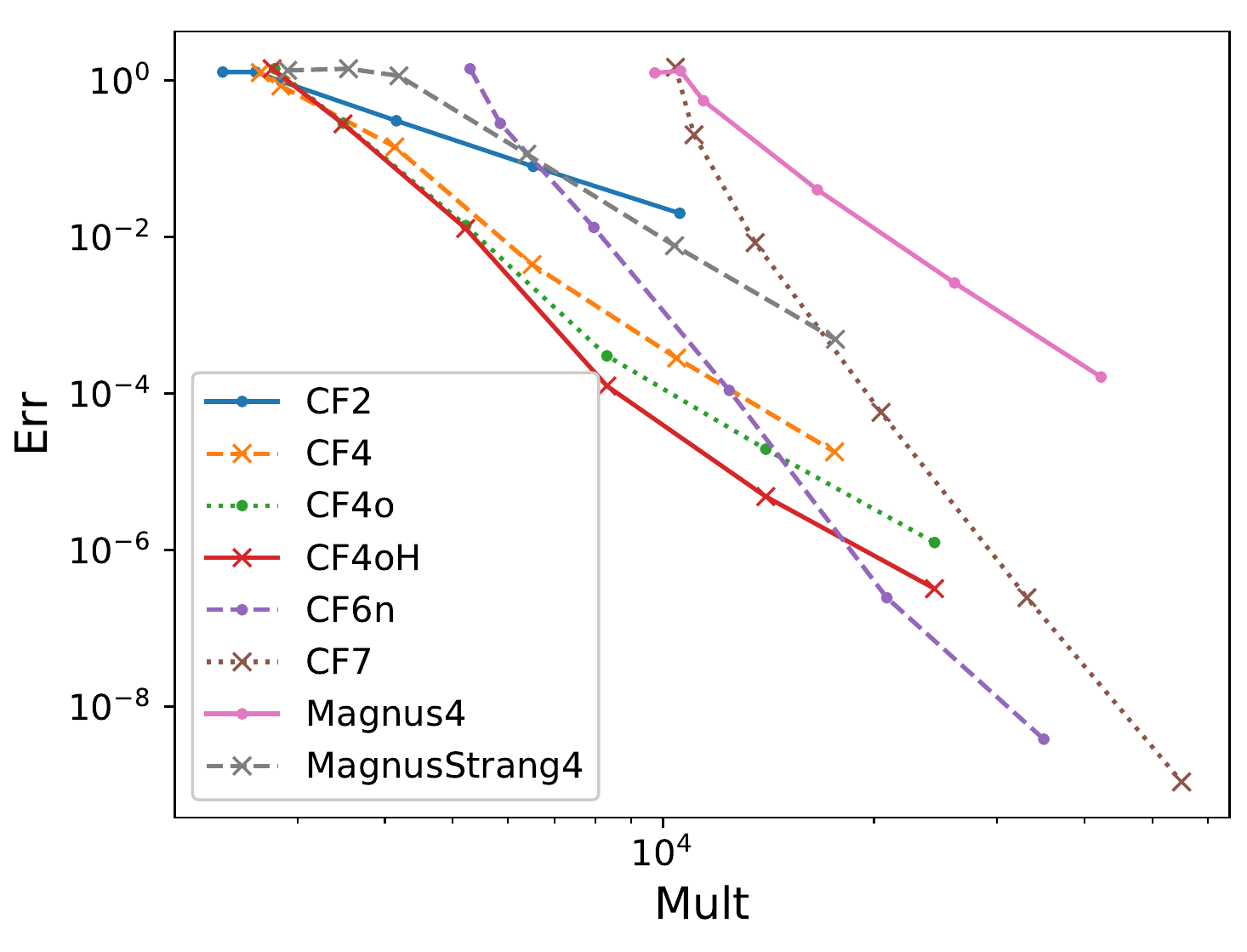} \quad \includegraphics[width=5.8cm]{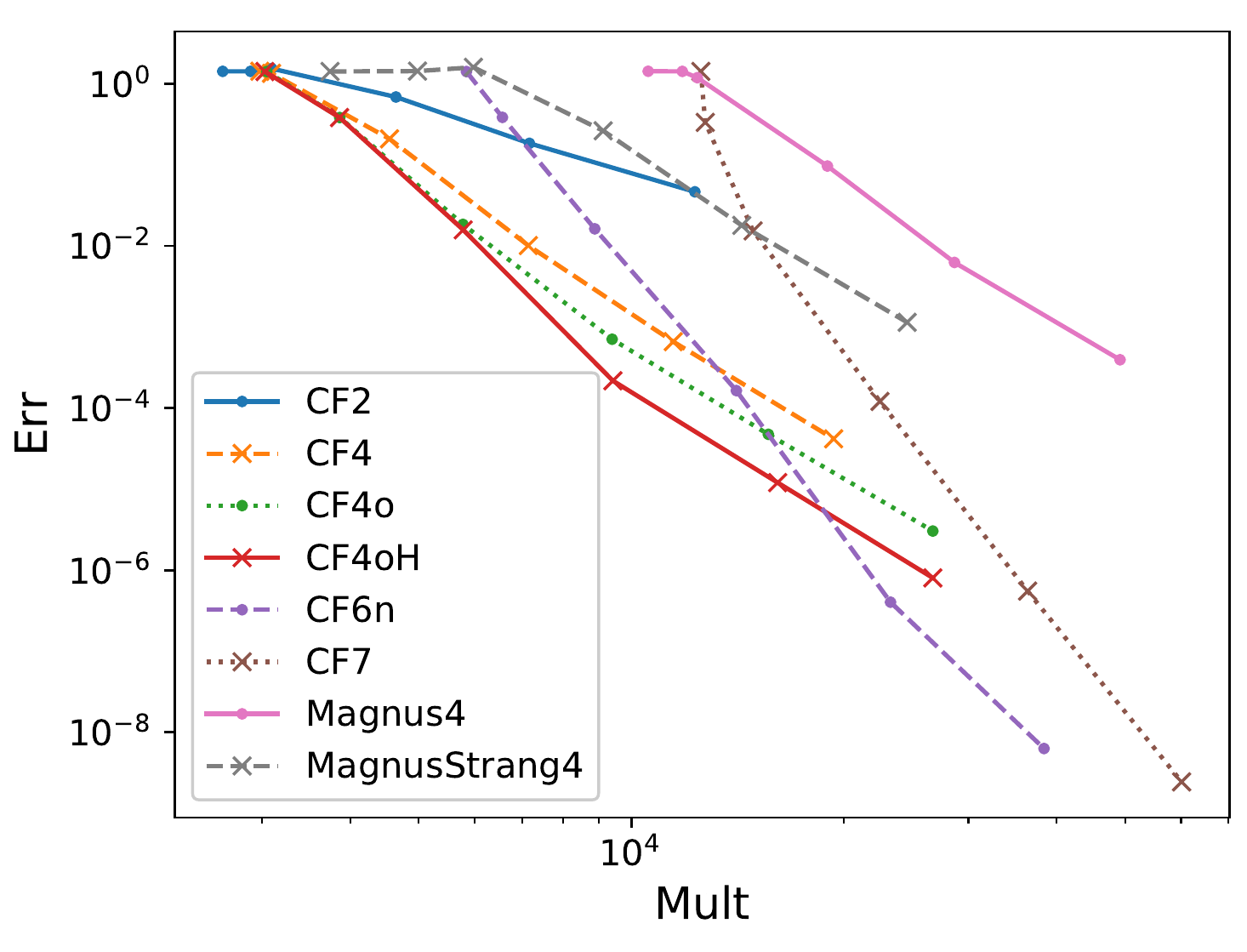} \\
\includegraphics[width=5.8cm]{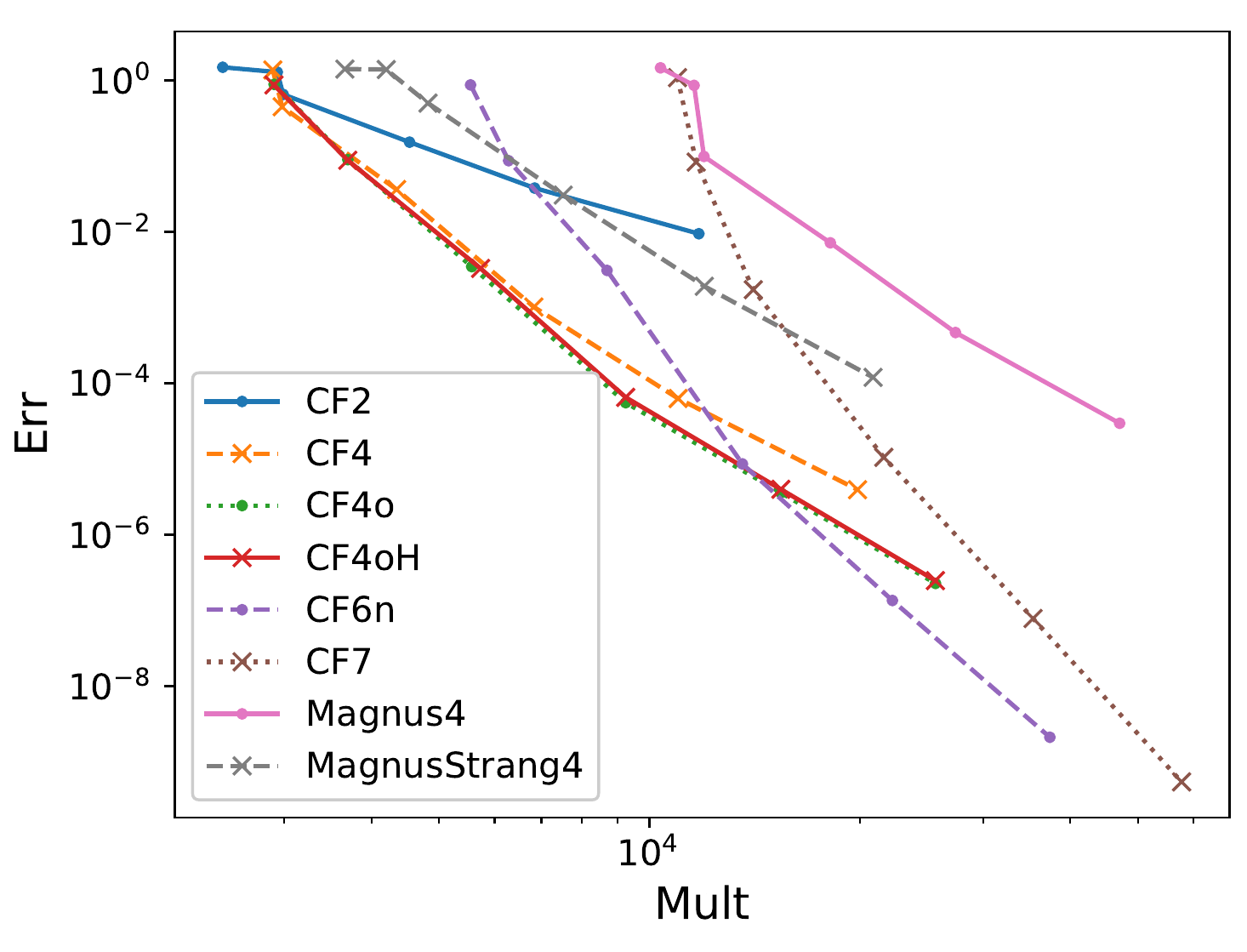} \quad \includegraphics[width=5.8cm]{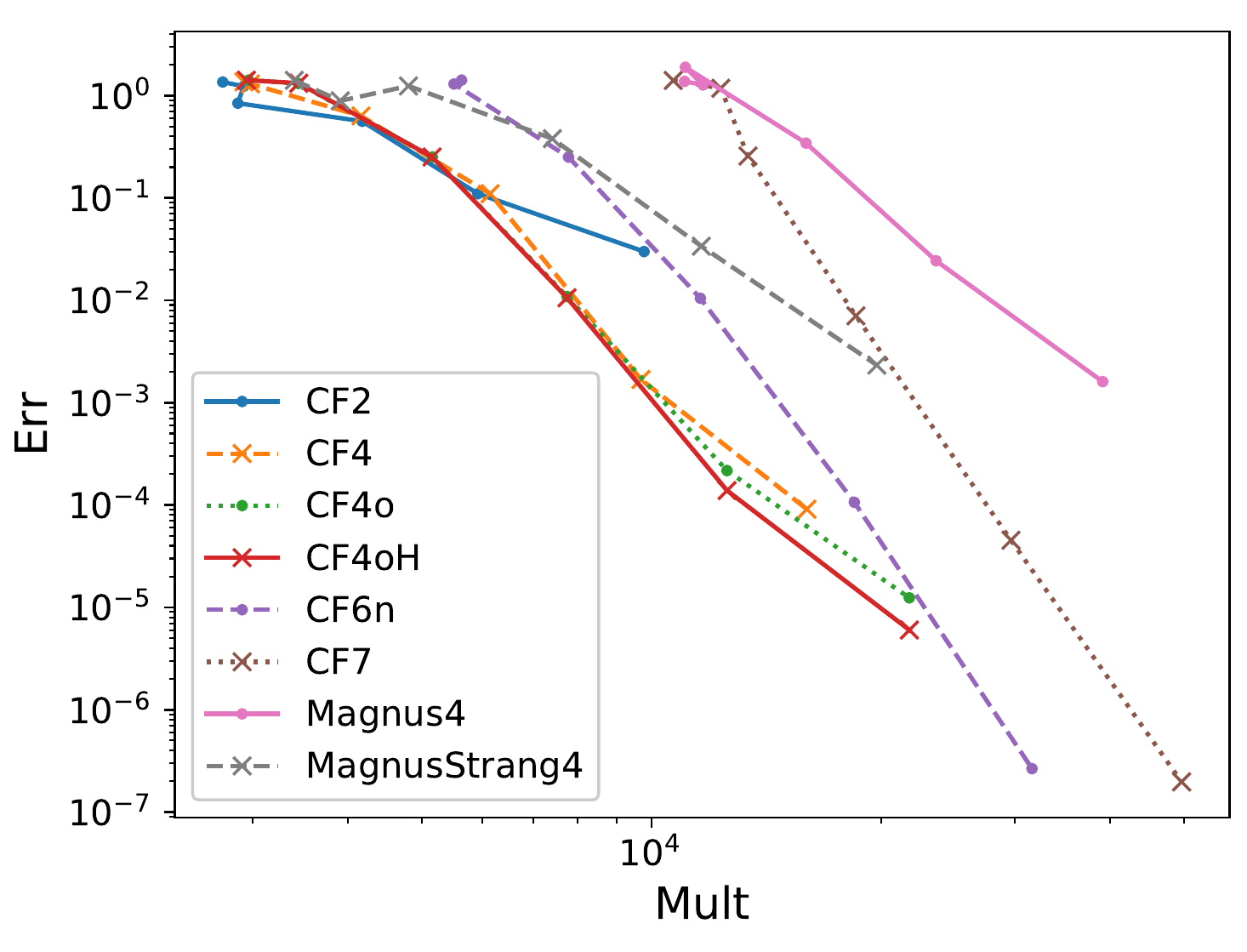} \\
\caption{$4\times3$ geometry, equidistant time-stepping. Error as a function of matrix--vector multiplications
for $\omega=3.5$ (top row), $\sigma_p=1$ (top left) and $\sigma_p=4$ (top right),
and for $\sigma_p=2$ (bottom row), $\omega=1.75$ (bottom left) and $\omega=7$ (bottom right).
\label{4x3variationequidistmult}}
\end{center}
\end{figure}

Figures~\ref{8x1variationequidist} and~\ref{8x1variationequidistmult} show the error
for equidistant time-stepping for the $2\times4$ geometry. Again, we observe an advantage for the highest-order
methods when high accuracy is sought, \texttt{CF4oH} is, as before, the most accurate fourth-order method.
Figure~\ref{4x3variationequidistmult} shows the same picture for the $4\times3$ geometry.

\begin{figure}[ht!]
\begin{center}
\includegraphics[width=5.8cm]{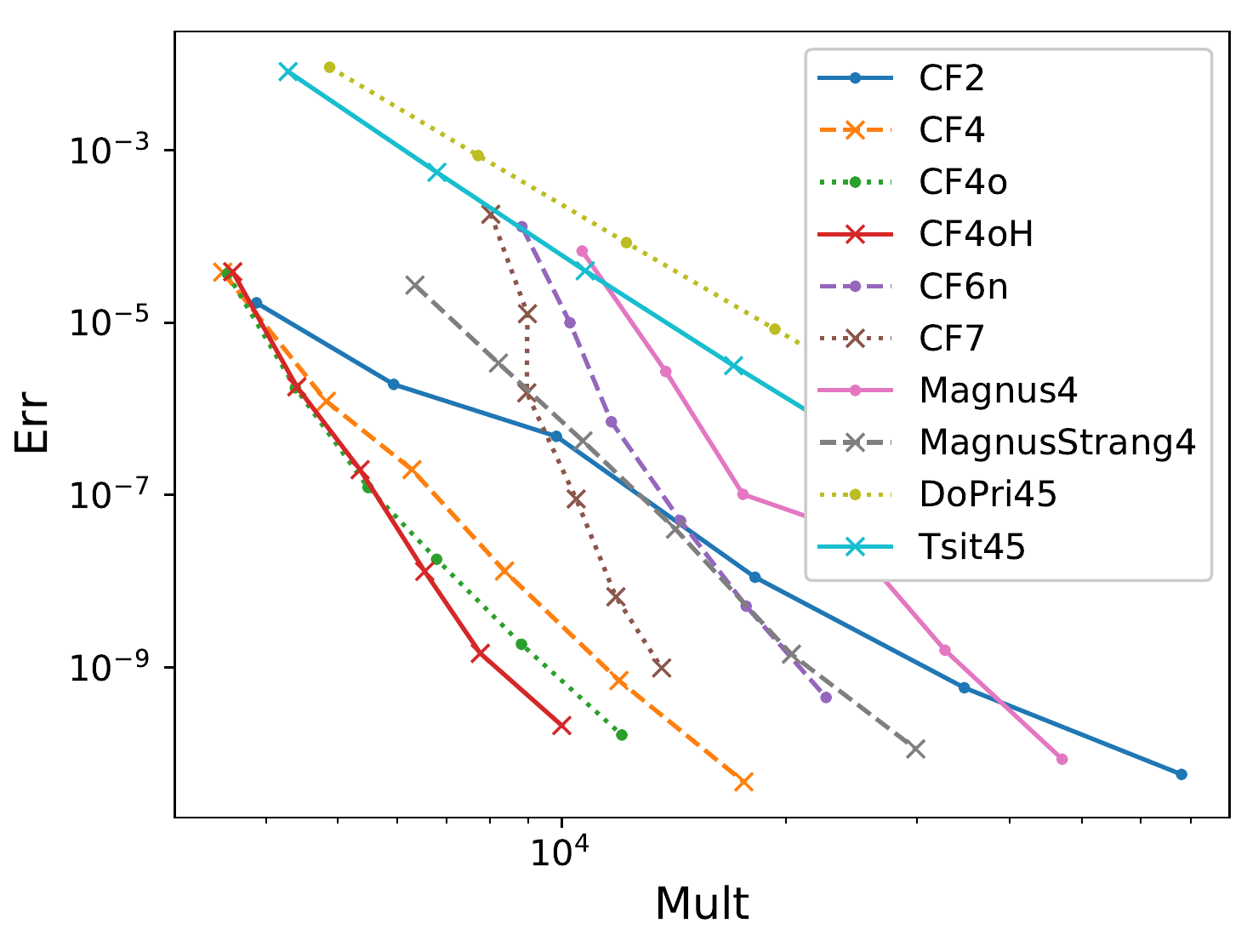} \quad \includegraphics[width=5.8cm]{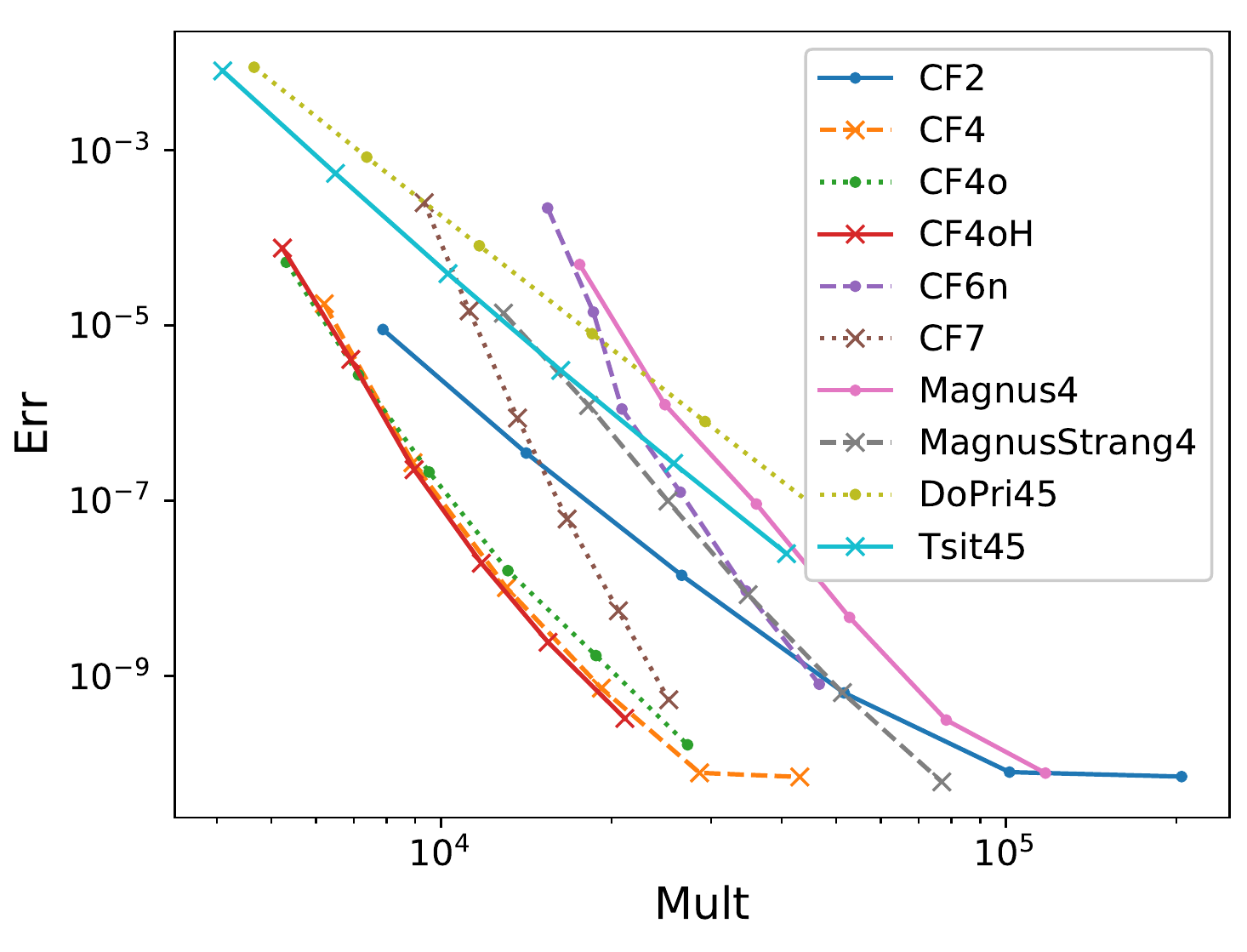} \\
\includegraphics[width=5.8cm]{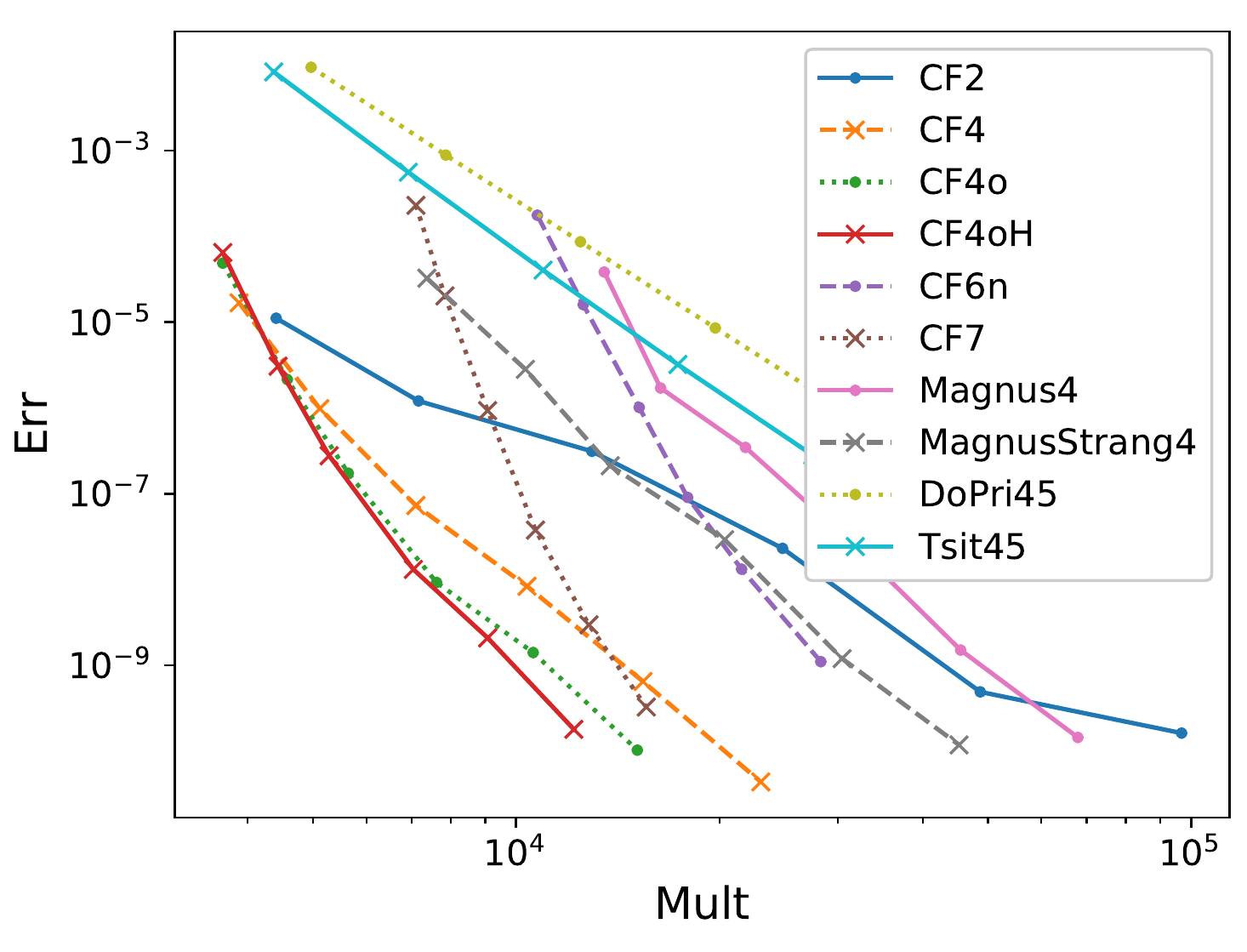} \quad \includegraphics[width=5.8cm]{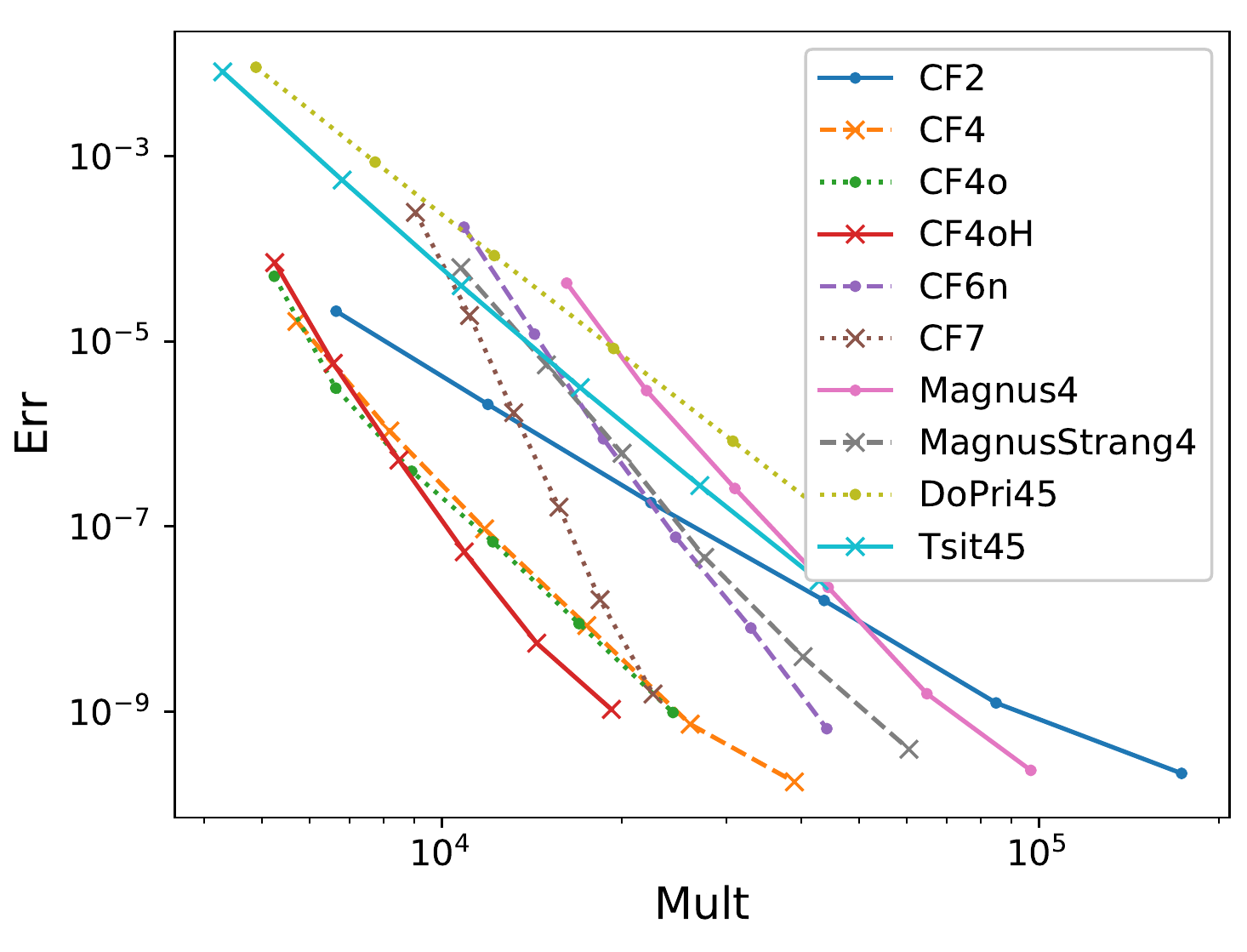} \\
\caption{$2\times4$ geometry, adaptive time-stepping. Error as a function of matrix-vector
multiplications for $\omega=3.5$ (top row), $\sigma_p=1$ (top left) and $\sigma_p=4$ (top right),
and for $\sigma_p=2$ (bottom row), $\omega=1.75$ (bottom left) and $\omega=7$ (bottom right).
\label{adapt8x1variation}}
\end{center}
\end{figure}

Finally, Figure~\ref{adapt8x1variation} shows the accuracy as a function of matrix--vector
multiplications for adaptive time-stepping for the $2\times4$ geometry. In this respect, \texttt{CF4oH}
is the most efficient choice. The advantage of adaptivity is quite pronounced for
these choices of parameters, a comparison of Figures~\ref{8x1variationequidistmult}
and~\ref{adapt8x1variation} shows that for a given number of matrix--vector multiplications,
the achieved accuracy is significantly higher in the adaptive integration.

\end{appendix}


\end{document}